\documentclass[11pt,a4paper]{article}
\usepackage{theorem,enumerate}
\usepackage{amsmath,latexsym,amssymb,amsfonts}
\usepackage{eucal}
\usepackage{color}
\usepackage{comment}
\theorembodyfont{\normalfont\slshape}
\newtheorem{thm}{Theorem}[section]
\newtheorem{cor}[thm]{Corollary}
\newtheorem{prop}[thm]{Proposition}
\newtheorem{lem}[thm]{Lemma}

\newtheorem{Ex}[thm]{Example}
\newtheorem{rem}[thm]{Remark}

\newtheorem{claim}{Claim}[section]

\newtheorem{Thm}{Theorem}

\newtheorem{conj}{Conjecture}

\newcommand{\proof}{\medbreak\noindent\textit{Proof.}\quad}

\newcommand{\qed}{{$\quad\square$\vs{3.6}}}

\newcommand{\vs}[1]{\vspace*{#1 mm}}

\def\AA{{ \mathcal{A}}}

\def\GG{{ \mathcal{G}}}

\def\II{{ \mathcal{I}}}
\def\JJ{{ \mathcal{J}}}

\def\LL{{ \mathcal{L}}}

\bfseries\normalfont

\title{An algebraic reduction of Hedetniemi's conjecture}
\author{
Ryoya Fukasaku\footnote{Faculty of Mathematics, Kyushu University, Motooka 744, Nishi-ku, Fukuoka 819-0395, Japan. \texttt{email:fukasaku@math.kyushu-u.ac.jp}}\and \
Michitaka Furuya\footnote{College of Liberal Arts and Sciences, Kitasato University, 1-15-1 Kitasato, Minami-ku, Sagamihara, Kanagawa 252-0373, Japan. \texttt{e-mail:michitaka.furuya@gmail.com}}\and \
Akihiro Higashitani\footnote{Department of Pure and Applied Mathematics, Graduate School of Information Science and Technology, Osaka University, Osaka 565-0871, Japan. \texttt{email:higashitani@ist.osaka-u.ac.jp}} 
}

\date{}

\begin{document}

\maketitle

\begin{abstract}
For a graph $G$, let $\chi (G)$ denote the chromatic number. 
In graph theory, the following famous conjecture posed by Hedetniemi has been studied:
For two graphs $G$ and $H$, $\chi (G\times H)=\min\{\chi (G),\chi (H)\}$, where $G \times H$ is the tensor product of $G$ and $H$. 
In this paper, we give a reduction of Hedetniemi's conjecture to an inclusion relation problem on ideals of polynomial rings, and we demonstrate computational experiments for partial solutions of Hedetniemi's conjecture along such a strategy using Gr\"{o}bner basis.
\end{abstract}



%%%%%%%%%%%%%%%%%%%%%%%%%%%%%%%%%%%%%%%%%%%%%%%%%%%%%%%%%%%%%%%%%%%%%%%%%%%%%%%%%%%%%%%%%%%%%%%%%%%%%%%%%%%%%%%%%%%%%%%%
%%%%%%%%%%%%%%%%%%%%%%%%%%%%%%%%%%%%%%%%%%%%%%%%%%%%%%%%%%%%%%%%%%%%%%%%%%%%%%%%%%%%%%%%%%%%%%%%%%%%%%%%%%%%%%%%%%%%%%%%
%%%%%%%%%%%%%%%%%%%%%%%%%%%%%%%%%%%%%%%%%%%%%%%%%%%%%%%%%%%%%%%%%%%%%%%%%%%%%%%%%%%%%%%%%%%%%%%%%%%%%%%%%%%%%%%%%%%%%%%%
\section{Introduction}\label{sec-intro}
%%%%%%%%%%%%%%%%%%%%%%%%%%%%%%%%%%%%%%%%%%%%%%%%%%%%%%%%%%%%%%%%%%%%%%%%%%%%%%%%%%%%%%%%%%%%%%%%%%%%%%%%%%%%%%%%%%%%%%%%
%%%%%%%%%%%%%%%%%%%%%%%%%%%%%%%%%%%%%%%%%%%%%%%%%%%%%%%%%%%%%%%%%%%%%%%%%%%%%%%%%%%%%%%%%%%%%%%%%%%%%%%%%%%%%%%%%%%%%%%%
%%%%%%%%%%%%%%%%%%%%%%%%%%%%%%%%%%%%%%%%%%%%%%%%%%%%%%%%%%%%%%%%%%%%%%%%%%%%%%%%%%%%%%%%%%%%%%%%%%%%%%%%%%%%%%%%%%%%%%%%

In this paper, we consider only finite undirected simple graphs.
Let $G$ be a graph.
Let $V(G)$ and $E(G)$ denote the {\it vertex set} and the {\it edge set} of $G$, respectively.
For $u\in V(G)$, let $N_{G}(u)$ and $d_{G}(u)$ denote the {\it neighborhood} and the {\it degree} of $u$, respectively; thus $N_{G}(u)=\{v\in V(G):uv\in E(G)\}$ and $d_{G}(u)=|N_{G}(u)|$.
Let $\delta (G)$ denote the {\it minimum degree} of $G$.
For $u\in V(G)$ and $X\subseteq V(G)$, let $\mbox{dist}_{G}(u,X)$ denote the length of a shortest path of $G$ joining $u$ and some vertex in $X$.
For $X\subseteq V(G)$, let $G[X]$ denote the subgraph of $G$ induced by $X$.
A $k$-subset $X\subseteq V(G)$ is called a {\it $k$-clique} (or just a {\it clique}) of $G$ if $G[X]$ is a complete graph.
Let $\omega (G)$ denote the largest positive integer $k$ such that $G$ contains a $k$-clique.
Let $K_{n}$ and $C_{n}$ denote the {\it complete graph} and the {\it cycle} of order $n$, respectively.
For a positive integer $k$, a mapping $c:V(G)\rightarrow [k]$ is a {\it proper $k$-coloring} of $G$ if $c(u)\neq c(v)$ for all adjacent vertices $u$ and $v$ of $G$, where $[k]=\{1,2,\ldots ,k\}$.
The smallest positive integer $k$ such that $G$ has a proper $k$-coloring is called the {\it chromatic number} of $G$, and it is denoted by $\chi (G)$.
For terms and symbols not defined here, we refer the reader to \cite{BM}.

The product operations of graphs have been widely studied because they can produce important illustrations for many graph properties.
The readers might find many interesting results in, for example, \cite{IK}.
In the deep studies for products, some primitive (but essential) problems and conjectures were posed.
In this paper, we focus on a classical conjecture concerning the chromatic number of a product of graphs.
Let $G$ and $H$ be two graphs.
The {\it tensor product} $G\times H$ of $G$ and $H$ is the graph on $V(G)\times V(H)$ such that two vertices $(u,v)$ and $(u',v')$ are adjacent in $G\times H$ if and only if $uu'\in E(G)$ and $vv'\in E(H)$.
For a proper $k$-coloring $c$ of $G$, the mapping $c_{0}:V(G)\times V(H)\rightarrow [k]$ with $c_{0}(u,v)=c(u)~(u\in V(G),~v\in V(H))$ is clearly a proper $k$-coloring of $G\times H$.
By the symmetry of $G$ and $H$, this leads to
\begin{align}
\chi (G\times H)\leq \min\{\chi (G),\chi (H)\}.\label{eq-intro-1}
\end{align}
Hedetniemi~\cite{H} conjectured that the equality in (\ref{eq-intro-1}) always holds for all graphs $G$ and $H$ as follows.

\begin{conj}[Hedetniemi~\cite{H}]%%%%%%%%%%%%%%%%%%%%%%%%%%%%%%%%%%%%%%%%%%%%%%%%%%%%%%%%%%%%%%%%%%%%%%%%%%%%%%%%%%%%%%%
\label{conj1}
Let $G$ and $H$ be graphs.
Then $\chi (G\times H)=\min\{\chi (G),\chi (H)\}$.
\end{conj}
%%%%%%%%%%%%%%%%%%%%%%%%%%%%%%%%%%%%%%%%%%%%%%%%%%%%%%%%%%%%%%%%%%%%%%%%%%%%%%%%%%%%%%%%%%%%%%%%%%%%%%%%%%%%%%%%%%%%%%%%

Conjecture~\ref{conj1} has been studied for more than 50 years, and some approaches for the conjecture and its analogies were also studied (see surveys \cite{T,Z}).
Hedetniemi~\cite{H} verified that Conjecture~\ref{conj1} is true for the case where $\min\{\chi (G),\chi (H)\}\leq 3$, and El-Zahar and Sauer~\cite{ES} proved that Conjecture~\ref{conj1} is true if $\min\{\chi (G),\chi (H)\}=4$.
On the other hand, Shitov~\cite{Sh} recently constructed counterexamples for Conjecture~\ref{conj1}, and he proved that if an integer $k$ is sufficiently large, then there exist infinitely many pairs $(G,H)$ of graphs such that $\min\{\chi (G),\chi (H)\}>k$ and $\chi (G\times H)=k$.
We remark that Conjecture~\ref{conj1} is still open for the case where $\min\{\chi (G),\chi (H)\}$ is small. 
Our aim of this paper is to propose a new effective approach to solve Conjecture~\ref{conj1} for small $k$'s; a reduction of the conjecture to an inclusion relation problem on ideals of polynomial rings.
Indeed, we demonstrate computational experiments for partial solutions of the conjecture using Gr\"{o}bner basis.

We are inspired from the results given by Margulies and Hicks~\cite{MH} concerning Vizing's conjecture, that is a conjecture on the domination number of Cartesian product of graphs.
They also reduced Vizing's conjecture to an inclusion relation problem of ideals.
However, the chromatic number and the domination number have major difference for the criticality.
When we consider a reduction of graph-theoritical problems to an inclusion relation of ideal, the criticality concerning edge-deletion or vertex-deletion is a useful tool. 
For a given graph $G$, although a subgraph of $G$ might have larger domination number than $G$, deleting a vertex or an edge cannot increase the chromatic number, that is, every graph contains a subgraph with a criticality for the chromatic number.
This property gives a strong advantage if we adopt the reduction strategy to Hedetniemi's conjecture.

This paper is organized as follows:
In Section~\ref{sec-prel}, we list some known results concerning Conjecture~\ref{conj1}.
In order to clear the standpoint of the cases treated in our computational experiments, we indicate the cases which force Conjecture~\ref{conj1} to be true from known results in Section~\ref{sec-small}.
Some results proved in Section~\ref{sec-small} might been known, but to keep the paper self-contained we give their proofs. 
Thus readers not interested in its detail are advised to skip the proof.
The main results are in Section~\ref{sec-main}.
In Subsection~\ref{sec-main1}, we reduce Conjecture~\ref{conj1} to a problem concerning graphs with the criticality for chromatic numbers. 
Using the reduction, we further reduce the conjecture to an inclusion relation problem on ideals of polynomial rings in Subsection~\ref{sec-main2}.
In Subsection~\ref{sec-main3}, more feasible reductions for computer analysis are considered.
In Section~\ref{sec-nume}, we give computational experiments along the strategy developed in Section~\ref{sec-main}.

%%%%%%%%%%%%%%%%%%%%%%%%%%%%%%%%%%%%%%%%%%%%%%%%%%%%%%%%%%%%%%%%%%%%%%%%%%%%%%%%%%%%%%%%%%%%%%%%%%%%%%%%%%%%%%%%%%%%%%%%
%%%%%%%%%%%%%%%%%%%%%%%%%%%%%%%%%%%%%%%%%%%%%%%%%%%%%%%%%%%%%%%%%%%%%%%%%%%%%%%%%%%%%%%%%%%%%%%%%%%%%%%%%%%%%%%%%%%%%%%%
%%%%%%%%%%%%%%%%%%%%%%%%%%%%%%%%%%%%%%%%%%%%%%%%%%%%%%%%%%%%%%%%%%%%%%%%%%%%%%%%%%%%%%%%%%%%%%%%%%%%%%%%%%%%%%%%%%%%%%%%
\section{Preliminary results}\label{sec-prel}
%%%%%%%%%%%%%%%%%%%%%%%%%%%%%%%%%%%%%%%%%%%%%%%%%%%%%%%%%%%%%%%%%%%%%%%%%%%%%%%%%%%%%%%%%%%%%%%%%%%%%%%%%%%%%%%%%%%%%%%%
%%%%%%%%%%%%%%%%%%%%%%%%%%%%%%%%%%%%%%%%%%%%%%%%%%%%%%%%%%%%%%%%%%%%%%%%%%%%%%%%%%%%%%%%%%%%%%%%%%%%%%%%%%%%%%%%%%%%%%%%
%%%%%%%%%%%%%%%%%%%%%%%%%%%%%%%%%%%%%%%%%%%%%%%%%%%%%%%%%%%%%%%%%%%%%%%%%%%%%%%%%%%%%%%%%%%%%%%%%%%%%%%%%%%%%%%%%%%%%%%%

In this section, we list some useful results for our argument.

\begin{Thm}[Burr, Erd\H{o}s and Lov\'{a}sz~\cite{BEL}]%%%%%%%%%%%%%%%%%%%%%%%%%%%%%%%%%%%%%%%%%%%%%%%%%%%%%%%%%%%%%%%%%%
\label{Thm-prel-B}
Let $k\geq 2$ be an integer.
Let $G$ and $H$ be graphs with $\chi (G)=\chi (H)=k$, and suppose that each vertex of $G$ belongs to a $(k-1)$-clique of $G$.
Then $\chi (G\times H)=k$.
\end{Thm}
%%%%%%%%%%%%%%%%%%%%%%%%%%%%%%%%%%%%%%%%%%%%%%%%%%%%%%%%%%%%%%%%%%%%%%%%%%%%%%%%%%%%%%%%%%%%%%%%%%%%%%%%%%%%%%%%%%%%%%%%

\begin{Thm}[Duffus, Sands and Woodrow~\cite{DSW}; Welzl~\cite{W}]%%%%%%%%%%%%%%%%%%%%%%%%%%%%%%%%%%%%%%%%%%%%%%%%%%%%%%%
\label{Thm-prel-C}
Let $k\geq 2$ be an integer.
Let $G$ and $H$ be graphs with $\chi (G)=\chi (H)=k$, and suppose that both $G$ and $H$ contain $(k-1)$-cliques.
Then $\chi (G\times H)=k$.
\end{Thm}
%%%%%%%%%%%%%%%%%%%%%%%%%%%%%%%%%%%%%%%%%%%%%%%%%%%%%%%%%%%%%%%%%%%%%%%%%%%%%%%%%%%%%%%%%%%%%%%%%%%%%%%%%%%%%%%%%%%%%%%%

As we mentioned in Section~\ref{sec-intro}, the criticality for chromatic number plays a crucial role in this paper.
Thus we next focus on such a concept and related results.

A graph $G$ is said to be {\it $k$-critical} if $\chi (G)=k$ and $\chi (G')\leq k-1$ for all subgraphs $G'$ of $G$ with $G'\neq G$.
In many papers, edge-critical graphs (i.e., graphs $G$ with $\chi (G-e)\leq \chi (G)-1$ for all $e\in E(G)$) and vertex-critical graphs (i.e., graphs $G$ with $\chi (G-u)\leq \chi (G)-1$ for all $u\in V(G)$) are individually considered.
Note that the concept of critical graphs defined above contains such two criticality concepts.
It is clear that $K_{k}$ is the unique $k$-critical graph if $k\in \{1,2\}$.
Furthermore, a graph is $3$-critical if and only if the graph is an odd cycle.
On the other hand, nobody knows an explicit characterization of $4$-critical graphs, and $4$-critical graphs have been studied. 

For two vertex-disjoint graphs $G_{1}$ and $G_{2}$, the {\it join} of $G_{1}$ and $G_{2}$, denoted by $G_{1}+G_{2}$, is obtained from $G_{1}$ and $G_{2}$ by joining each vertex of $G_{1}$ to all vertices of $G_{2}$.
A graph $G$ is {\it decomposable} if the complement $\overline{G}$ of $G$ is disconnected.
A non-decomposable graph is said to be {\it indecomposable}.
We can easily verify that a $k$-critical graph $G$ is decomposable if and only if $G$ is the join of a $k_{1}$-critical graph $G_{1}$ and a $k_{2}$-critical graph $G_{2}$ with $k_{1}+k_{2}=k$.
On the other hand, indecomposable critical graphs have many vertices as follows (here the second statement was proved by Gallai~\cite{G}):

\begin{Thm}[Stehl\'{i}k~\cite{S}]%%%%%%%%%%%%%%%%%%%%%%%%%%%%%%%%%%%%%%%%%%%%%%%%%%%%%%%%%%%%%%%%%%%%%%%%%%%%%%%%%%%%%%%
\label{Thm-prel-A}
Let $k\geq 3$ be an integer, and let $G$ be an indecomposable $k$-critical graph.
Then for any $u\in V(G)$, $G-u$ has a proper $(k-1)$-coloring such that every color class contains at least two vertices.
In particular, $|V(G)|\geq 2k-1$.
\end{Thm}
%%%%%%%%%%%%%%%%%%%%%%%%%%%%%%%%%%%%%%%%%%%%%%%%%%%%%%%%%%%%%%%%%%%%%%%%%%%%%%%%%%%%%%%%%%%%%%%%%%%%%%%%%%%%%%%%%%%%%%%%

Furthermore, the following result closely related to the $k$-criticality is well-known and we can find it in many textbooks of graph theory (for example, in \cite[Theorem~14.7]{BM}).

\begin{Thm}%%%%%%%%%%%%%%%%%%%%%%%%%%%%%%%%%%%%%%%%%%%%%%%%%%%%%%%%%%%%%%%%%%%%%%%%%%%%%%%%%%%%%%%%%%%%%%%%%%%%%%%%%%%%%
\label{prop-verc-mindeg}
For a positive integer $k$, every $k$-critical graph $G$ satisfies $\delta (G)\geq k-1$.
\end{Thm}
%%%%%%%%%%%%%%%%%%%%%%%%%%%%%%%%%%%%%%%%%%%%%%%%%%%%%%%%%%%%%%%%%%%%%%%%%%%%%%%%%%%%%%%%%%%%%%%%%%%%%%%%%%%%%%%%%%%%%%%%

%%%%%%%%%%%%%%%%%%%%%%%%%%%%%%%%%%%%%%%%%%%%%%%%%%%%%%%%%%%%%%%%%%%%%%%%%%%%%%%%%%%%%%%%%%%%%%%%%%%%%%%%%%%%%%%%%%%%%%%%
%%%%%%%%%%%%%%%%%%%%%%%%%%%%%%%%%%%%%%%%%%%%%%%%%%%%%%%%%%%%%%%%%%%%%%%%%%%%%%%%%%%%%%%%%%%%%%%%%%%%%%%%%%%%%%%%%%%%%%%%
%%%%%%%%%%%%%%%%%%%%%%%%%%%%%%%%%%%%%%%%%%%%%%%%%%%%%%%%%%%%%%%%%%%%%%%%%%%%%%%%%%%%%%%%%%%%%%%%%%%%%%%%%%%%%%%%%%%%%%%%
\section{Hedetniemi's conjecture for small graphs}\label{sec-small}
%%%%%%%%%%%%%%%%%%%%%%%%%%%%%%%%%%%%%%%%%%%%%%%%%%%%%%%%%%%%%%%%%%%%%%%%%%%%%%%%%%%%%%%%%%%%%%%%%%%%%%%%%%%%%%%%%%%%%%%%
%%%%%%%%%%%%%%%%%%%%%%%%%%%%%%%%%%%%%%%%%%%%%%%%%%%%%%%%%%%%%%%%%%%%%%%%%%%%%%%%%%%%%%%%%%%%%%%%%%%%%%%%%%%%%%%%%%%%%%%%
%%%%%%%%%%%%%%%%%%%%%%%%%%%%%%%%%%%%%%%%%%%%%%%%%%%%%%%%%%%%%%%%%%%%%%%%%%%%%%%%%%%%%%%%%%%%%%%%%%%%%%%%%%%%%%%%%%%%%%%%

In this section, we focus on small graphs $G$ and $H$ satisfying Conjecture~\ref{conj1} and finally prove the following theorem.

\begin{thm}%%%%%%%%%%%%%%%%%%%%%%%%%%%%%%%%%%%%%%%%%%%%%%%%%%%%%%%%%%%%%%%%%%%%%%%%%%%%%%%%%%%%%%%%%%%%%%%%%%%%%%%%%%%%%
\label{thm-small-main}
Let $k\geq 5$ be an integer.
Let $G$ and $H$ be graphs with $\min\{\chi (G),\chi (H)\}=k$, and suppose that
\begin{enumerate}[{\upshape(i)}]
\item
$\min\{|V(G)|,|V(H)|\}\leq k+2$;
\item
$|V(G)|=|V(H)|=k+3$; or
\item
$k=5$ and $(|V(G)|,|V(H)|)\in \{(8,8),(8,9),(8,10),(9,8),(10,8)\}$.
\end{enumerate}
Then $\chi (G\times H)=k$.
\end{thm}
%%%%%%%%%%%%%%%%%%%%%%%%%%%%%%%%%%%%%%%%%%%%%%%%%%%%%%%%%%%%%%%%%%%%%%%%%%%%%%%%%%%%%%%%%%%%%%%%%%%%%%%%%%%%%%%%%%%%%%%%

By Theorem~\ref{thm-small-main}, the first nontrivial cases for Conjecture~\ref{conj1} are
\begin{enumerate}[{$\bullet $}]
\item
$|V(G)|=8$ and $|V(H)|=11$ if $k=5$; and
\item
$|V(G)|=k+3$ and $|V(H)|=k+4$ if $k\geq 6$.
\end{enumerate}
We can refine the latter case as follows.
(Here, for a graph $H$, we regard $K_{0}+H$ as $H$.)

\begin{figure}[h]
\begin{center}
\input{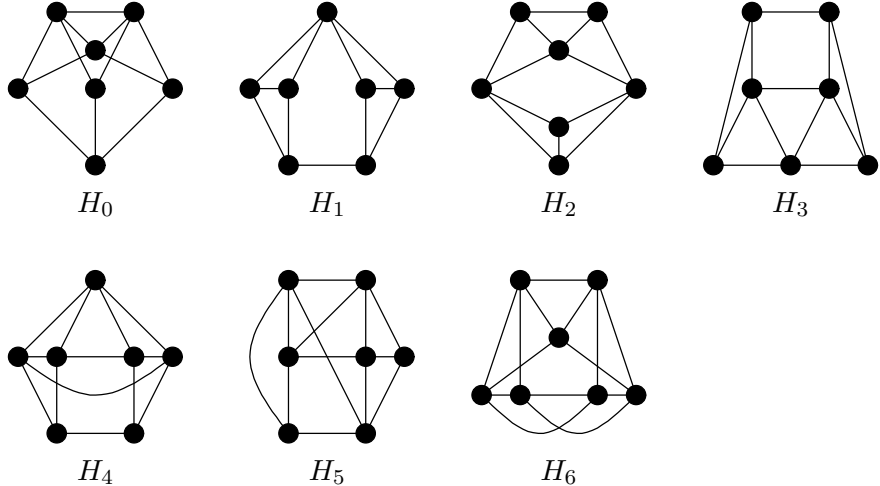}
\caption{The $4$-critical graphs of order $7$.}
\label{f0}
\end{center}
\end{figure}

\begin{thm}%%%%%%%%%%%%%%%%%%%%%%%%%%%%%%%%%%%%%%%%%%%%%%%%%%%%%%%%%%%%%%%%%%%%%%%%%%%%%%%%%%%%%%%%%%%%%%%%%%%%%%%%%%%%%
\label{thm-small-main2}
Let $k\geq 6$ be an integer.
Then all graphs $G$ and $H$ with $|V(G)|\leq k+3$, $|V(H)|\leq k+4$ and $\min\{\chi (G),\chi (H)\}=k$ satisfy $\chi (G\times H)=k$ if and only if 
$$\chi ((K_{k-4}+H_{0})\times (K_{k-6}+C_{5}+C_{5}))=k,$$
where $H_0$ denotes the graph depicted in Figure~\ref{f0}. 
\end{thm}
%%%%%%%%%%%%%%%%%%%%%%%%%%%%%%%%%%%%%%%%%%%%%%%%%%%%%%%%%%%%%%%%%%%%%%%%%%%%%%%%%%%%%%%%%%%%%%%%%%%%%%%%%%%%%%%%%%%%%%%%

The following theorem is a useful tool in the proof of our argument.

\begin{Thm}[Chv\'{a}tal~\cite{Ch}; Jensen and Royle~\cite{JR}]%%%%%%%%%%%%%%%%%%%%%%%%%%%%%%%%%%%%%%%%%%%%%%%%%%%%%%%%%%
\label{Thm-smallcolor}
For $k\in \{4,5\}$, if a $K_{k-1}$-free graph $G$ satisfies $\chi (G)=k$, then $|V(G)|\geq 11$.
\end{Thm}
%%%%%%%%%%%%%%%%%%%%%%%%%%%%%%%%%%%%%%%%%%%%%%%%%%%%%%%%%%%%%%%%%%%%%%%%%%%%%%%%%%%%%%%%%%%%%%%%%%%%%%%%%%%%%%%%%%%%%%%%

We first prove that $K_{k}$ is the unique $k$-critical graph of order at most $k+1$.

\begin{lem}%%%%%%%%%%%%%%%%%%%%%%%%%%%%%%%%%%%%%%%%%%%%%%%%%%%%%%%%%%%%%%%%%%%%%%%%%%%%%%%%%%%%%%%%%%%%%%%%%%%%%%%%%%%%%
\label{lem-small-00}
Let $k\geq 1$ be an integer, and let $G$ be a $k$-critical graph of order at most $k+1$.
Then $G=K_{k}$.
\end{lem}
%%%%%%%%%%%%%%%%%%%%%%%%%%%%%%%%%%%%%%%%%%%%%%%%%%%%%%%%%%%%%%%%%%%%%%%%%%%%%%%%%%%%%%%%%%%%%%%%%%%%%%%%%%%%%%%%%%%%%%%%
\proof
It is clear that if $|V(G)|\leq k$, then $G=K_{k}$.
Thus it suffices to show that $|V(G)|\neq k+1$.
By way of contradiction, suppose that $|V(G)|=k+1$.
Let $c$ be a proper $k$-coloring of $G$.
We may assume that $|c^{-1}(i)|=1$ for every $i~(1\leq i\leq k-1)$ (and so $|c^{-1}(k)|=2$).
Note that $G[\bigcup _{i=1}^{k-1}c^{-1}(i)]$ is a complete graph.
Write $c^{-1}(i)=\{u_{i}\}$ for each $i~(1\leq i\leq k-1)$ and $c^{-1}(k)=\{v_{1},v_{2}\}$.
If for each $j\in \{1,2\}$, there exists a vertex $w_{j}\in \{u_{i}:1\leq i\leq k-1\}$ with $v_{j}w_{j}\notin E(G)$, 
then the mapping $c':V(G)\rightarrow [k-1]$ with
$$
c'(a)=
\begin{cases}
c(a) & (a\notin \{v_{1},v_{2}\})\\
c(w_{j}) & (a=v_{j})
\end{cases}
$$
is a proper $(k-1)$-coloring of $G$, which contradicts the fact that $\chi (G)=k$.
Thus, without loss of generality, we may assume that $N_{G}(v_{1})=V(G)\setminus \{v_{1},v_{2}\}$.
Then $G-v_{2}$ is a complete graph of order $k$, and so $\chi (G-v_{2})=k$, which contradicts the fact that $G$ is $k$-critical.
\qed

Let $\AA_{4}$ be the family of $4$-critical graphs of order $7$.
Then every $k$-critical graph with at most $k+3$ vertices can be characterized as follows.

\begin{lem}%%%%%%%%%%%%%%%%%%%%%%%%%%%%%%%%%%%%%%%%%%%%%%%%%%%%%%%%%%%%%%%%%%%%%%%%%%%%%%%%%%%%%%%%%%%%%%%%%%%%%%%%%%%%%
\label{lem-small-1}
For an integer $k\geq 3$, a graph $G$ of order at most $k+3$ is $k$-critical if and only if
\begin{enumerate}[{\upshape(i)}]
\item
$G=K_{k}$;
\item
$G=K_{k-3}+C_{5}$; or
\item
$k\geq 4$ and $G=K_{k-4}+A$ for a graph $A\in \AA_{4}$.
\end{enumerate}
\end{lem}
%%%%%%%%%%%%%%%%%%%%%%%%%%%%%%%%%%%%%%%%%%%%%%%%%%%%%%%%%%%%%%%%%%%%%%%%%%%%%%%%%%%%%%%%%%%%%%%%%%%%%%%%%%%%%%%%%%%%%%%%
\proof
The ``if'' part is trivial. Thus we show the ``only if'' part.
Since a graph is $3$-critical if and only if it is an odd cycle, the lemma holds for $k=3$.
Thus we may assume that $k\geq 4$.

\begin{claim}%%%%%%%%%%%%%%%%%%%%%%%%%%%%%%%%%%%%%%%%%%%%%%%%%%%%%%%%%%%%%%%%%%%%%%%%%%%%%%%%%%%%%%%%%%%%%%%%%%%%%%%%%%%
\label{cl-lem-small-1}
Let $l\geq 4$ be an integer, and let $H$ be an $l$-critical graph of order at most $l+3$.
Then either $(l,|V(H)|)=(4,7)$ or $H=K_{l-l_{0}}+H'$ for an $l_{0}$-critical graph $H'$ with $1\leq l_{0}\leq l-1$.
\end{claim}
%%%%%%%%%%%%%%%%%%%%%%%%%%%%%%%%%%%%%%%%%%%%%%%%%%%%%%%%%%%%%%%%%%%%%%%%%%%%%%%%%%%%%%%%%%%%%%%%%%%%%%%%%%%%%%%%%%%%%%%%
\proof
Note that $|V(H)|\leq l+3$ and $l+3 \leq 2l-1$ (i.e. $l \geq 4$) holds by our assumption. 
Hence, we have $|V(H)|\leq l+3 \leq 2l-1$ and all the equalities hold if and only if $(l,|V(H)|)=(4,7)$. 
Thus, to prove the claim, we may assume $|V(H)|<2l-1$. 
Then by Theorem~\ref{Thm-prel-A}, $H$ is decomposable, and hence $H=H_{1}+H_{2}$ for an $l_{1}$-critical graph $H_{1}$ and $l_{2}$-critical graph $H_{2}$ with $l_{1}+l_{2}=l$.
If both $H_{1}$ and $H_{2}$ are non-complete, then it follows from Lemma~\ref{lem-small-00} that $|V(H_{i})|\geq l_{i}+2~(i\in \{1,2\})$, and so $|V(H)|=|V(H_{1})|+|V(H_{2})|\geq (l_{1}+2)+(l_{2}+2)>l+3$, which is a contradiction.
Thus we may assume that $H_{1}$ is complete.
Then $H_{1}=K_{l_{1}}=K_{l-l_{2}}$, as desired.
\qed

By Lemma~\ref{lem-small-00}, we may assume that $|V(G)|\in \{k+2,k+3\}$, and so $G$ is non-complete.
If $(k,|V(G)|)=(4,7)$, then (iii) holds.
Thus we may assume that $(k,|V(G)|)\neq (4,7)$.
Then by Claim~\ref{cl-lem-small-1}, $G=K_{k-k_{0}}+G'$ for a $k_{0}$-critical graph $G'$ with $1\leq k_{0}\leq k-1$.
Note that $G'$ is non-complete.
Choose $k_{0}$ and $G'$ so that $|V(G')|$ is as small as possible.

In the case $k_{0}\leq 3$, since $G'$ is non-complete, $k_{0}=3$ and $G'$ is an odd cycle of order at least $5$.
Since $k-3=k-k_{0}=|V(G)\setminus V(G')|\in \{k+2-|V(G')|,k+3-|V(G')|\}$, it follows that $|V(G)|=k+2$ and $|V(G')|=5$, i.e., $G=K_{k-3}+C_{5}$, which implies (ii).

Let $k_{0}\geq 4$. 
Note that $k\geq 5$.
Since $|V(G')|=|V(G)|-(k-k_{0})\leq k+3-(k-k_{0})=k_{0}+3$, it follows from Claim~\ref{cl-lem-small-1} that either $(k_{0},|V(G')|)=(4,7)$ or $G'=K_{k_{0}-k_{1}}+G''$ for a $k_{1}$-critical graph $G''$ with $1\leq k_{1}\leq k_{0}-1$.
If the latter holds, then $G$ is the join of a complete graph of order $k-k_{1}$ and $G''$, which contradicts the choice of $k_{0}$ and $G'$.
Thus $(k_{0},|V(G')|)=(4,7)$, i.e., $G$ is the join of $K_{k-k_{0}}~(=K_{k-4})$ and $G'$ belonging to $\AA_{4}$, which implies (iii).
\qed

We will use Lemma~\ref{lem-small-1} to prove Theorem~\ref{thm-small-main}.
We can verify that $\AA_{4}$ consists of graphs $H_{0},H_{1},\ldots ,H_{6}$ depicted in Figure~\ref{f0}, and so Lemma~\ref{lem-small-1} gives a complete characterization of small $k$-critical graphs.
However, the characterization of $\AA_{4}$ might be proved by tedious argument (or computer search), and so we omit the detail.
Indeed, in order to prove Theorem~\ref{thm-small-main}, it suffices to prove a more restricted characterization as follows.

\begin{lem}%%%%%%%%%%%%%%%%%%%%%%%%%%%%%%%%%%%%%%%%%%%%%%%%%%%%%%%%%%%%%%%%%%%%%%%%%%%%%%%%%%%%%%%%%%%%%%%%%%%%%%%%%%%%%
\label{lem-small-2}
If a graph $G\in \AA_{4}$ has a vertex $u$ belonging to no triangle, then $G=H_{0}$.
\end{lem}
%%%%%%%%%%%%%%%%%%%%%%%%%%%%%%%%%%%%%%%%%%%%%%%%%%%%%%%%%%%%%%%%%%%%%%%%%%%%%%%%%%%%%%%%%%%%%%%%%%%%%%%%%%%%%%%%%%%%%%%%
\proof
By Theorem~\ref{prop-verc-mindeg}, we have
\begin{align}
\delta (G)\geq 3.\label{cond-lem-app-1}
\end{align}
Suppose that $G$ contains no triangle. 
Since $G$ has no proper $2$-coloring, $G$ contains an odd cycle of order at least $5$.
This together with \eqref{cond-lem-app-1} implies that $G$ contains an induced odd cycle $C$ of order $5$ and $N_{G}(v)\setminus V(C)\neq \emptyset $ for all $v\in V(C)$.
Since $|V(G)\setminus V(C)|=2$, a vertex in $V(G)\setminus V(C)$ is adjacent to two consecutive vertices on $C$, and so $G$ contains a triangle, which contradicts the assumption that $G$ contains no triangle.

Thus, $G$ contains a triangle $T=v_{1}v_{2}v_{3}v_{1}$.
Choose $T$ so that $\mbox{dist}_{G}(u,V(T))$ is as large as possible.
By the definition of $u$, we have $\mbox{dist}_{G}(u,V(T))\geq 1$.

\begin{figure}[htbp]
\begin{center}
%WinTpicVersion4.32a
{\unitlength 0.1in%
\begin{picture}(13.0000,9.0000)(1.5000,-10.5000)%
% CIRCLE 2 0 0 0 Black Black  
% 4 200 600 200 650 200 650 200 650
% 
\special{sh 1.000}%
\special{ia 200 600 50 50 0.0000000 6.2831853}%
\special{pn 8}%
\special{ar 200 600 50 50 0.0000000 6.2831853}%
% CIRCLE 2 0 0 0 Black Black  
% 4 600 600 600 650 600 650 600 650
% 
\special{sh 1.000}%
\special{ia 600 600 50 50 0.0000000 6.2831853}%
\special{pn 8}%
\special{ar 600 600 50 50 0.0000000 6.2831853}%
% CIRCLE 2 0 0 0 Black Black  
% 4 1000 600 1000 650 1000 650 1000 650
% 
\special{sh 1.000}%
\special{ia 1000 600 50 50 0.0000000 6.2831853}%
\special{pn 8}%
\special{ar 1000 600 50 50 0.0000000 6.2831853}%
% CIRCLE 2 0 0 0 Black Black  
% 4 1400 600 1400 650 1400 650 1400 650
% 
\special{sh 1.000}%
\special{ia 1400 600 50 50 0.0000000 6.2831853}%
\special{pn 8}%
\special{ar 1400 600 50 50 0.0000000 6.2831853}%
% CIRCLE 2 0 0 0 Black Black  
% 4 800 1000 800 1050 800 1050 800 1050
% 
\special{sh 1.000}%
\special{ia 800 1000 50 50 0.0000000 6.2831853}%
\special{pn 8}%
\special{ar 800 1000 50 50 0.0000000 6.2831853}%
% CIRCLE 2 0 0 0 Black Black  
% 4 400 200 400 250 400 250 400 250
% 
\special{sh 1.000}%
\special{ia 400 200 50 50 0.0000000 6.2831853}%
\special{pn 8}%
\special{ar 400 200 50 50 0.0000000 6.2831853}%
% CIRCLE 2 0 0 0 Black Black  
% 4 1200 200 1200 250 1200 250 1200 250
% 
\special{sh 1.000}%
\special{ia 1200 200 50 50 0.0000000 6.2831853}%
\special{pn 8}%
\special{ar 1200 200 50 50 0.0000000 6.2831853}%
% LINE 2 0 3 0 Black Black  
% 18 1200 200 400 200 400 200 200 600 200 600 1200 200 1200 200 600 600 600 600 400 200 400 200 1000 600 1000 600 1200 200 1200 200 1400 600 1400 600 400 200
% 
\special{pn 8}%
\special{pa 1200 200}%
\special{pa 400 200}%
\special{fp}%
\special{pa 400 200}%
\special{pa 200 600}%
\special{fp}%
\special{pa 200 600}%
\special{pa 1200 200}%
\special{fp}%
\special{pa 1200 200}%
\special{pa 600 600}%
\special{fp}%
\special{pa 600 600}%
\special{pa 400 200}%
\special{fp}%
\special{pa 400 200}%
\special{pa 1000 600}%
\special{fp}%
\special{pa 1000 600}%
\special{pa 1200 200}%
\special{fp}%
\special{pa 1200 200}%
\special{pa 1400 600}%
\special{fp}%
\special{pa 1400 600}%
\special{pa 400 200}%
\special{fp}%
% LINE 2 0 3 0 Black Black  
% 8 200 600 800 1000 800 1000 600 600 1000 600 800 1000 800 1000 1400 600
% 
\special{pn 8}%
\special{pa 200 600}%
\special{pa 800 1000}%
\special{fp}%
\special{pa 800 1000}%
\special{pa 600 600}%
\special{fp}%
\special{pa 1000 600}%
\special{pa 800 1000}%
\special{fp}%
\special{pa 800 1000}%
\special{pa 1400 600}%
\special{fp}%
% STR 2 0 3 0 Black White  
% 4 870 920 870 1020 1 0 0 0
% $u$
\put(8.7000,-10.2000){\makebox(0,0)[lt]{$u$}}%
\end{picture}}%
\caption{A graph appearing in the proof of Lemma~\ref{lem-small-2}.}
\label{f1}
\end{center}
\end{figure}
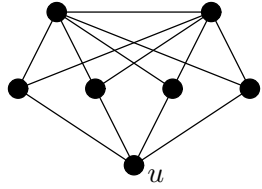

Suppose that $\mbox{dist}_{G}(u,V(T))=1$.
Note that $|N_{G}(u)\cap V(T)|=1$.
Without loss of generality, we may assume that $N_{G}(u)\cap V(T)=\{v_{1}\}$.
If $V(G)\setminus \{u,v_{2},v_{3}\}\subseteq N_{G}(u)$, then by \eqref{cond-lem-app-1} and the definition of $u$, 
$G$ is a graph depicted in Figure~\ref{f1}, and so $G$ has a proper $3$-coloring, 
which is a contradiction. Thus $V(G)\setminus \{u,v_{2},v_{3}\}\not\subseteq N_{G}(u)$.
This together with \eqref{cond-lem-app-1} implies that $d_{G}(u)=3$ and $V(G)\setminus (\{u,v_{1},v_{2}\}\cup N_{G}(u))$ contains exactly one vertex, say $x$.
%Write $N_{G}(u)=\{v_{1},w_{1},w_{2}\}$.
If $\{v_{2},v_{3}\}\subseteq N_{G}(x)$, then $\mbox{dist}_{G}(u,xv_{2}v_{3}x \})\geq 2$, which contradicts the choice of $T$. 
Thus, without loss of generality, we may assume that $xv_{2}\notin E(G)$.
Now we consider the mapping $c:V(G)\rightarrow [3]$ with
$$
c(a)=
\begin{cases}
1 & (a\in \{u,v_{2},x\})\\
2 & (a=v_{3})\\
3 & (a\in N_{G}(u)).
\end{cases}
$$
Then $c$ is a proper $3$-coloring of $G$ because $\{u,v_{2},x\}$ and $N_{G}(u)$ are independent sets of $G$, which is a contradiction.
Thus we may assume that $\mbox{dist}_{G}(u,V(T))\geq 2$.

Note that $V(G)=\{u\}\cup N_{G}(u)\cup V(T)$ and $d_{G}(u)=3$.
Write $N_{G}(u)=\{w_{1},w_{2},w_{3}\}$.
For $i\in \{1,2,3\}$, if $V(T)\subseteq N_{G}(w_{i})$, then the subgraph of $G$ induced by $V(T)\cup \{w_{i}\}$ is a complete graph of order $4$, and so $\chi (G-u)\geq 4$, which contradicts the $4$-criticality of $G$.
Since $\delta (G)\geq 3$, this implies that $|V(T)\cap N_{G}(w_{i})|=2$ for each $i\in \{1,2,3\}$.
Suppose that $V(T)\cap N_{G}(w_{i})=V(T)\cap N_{G}(w_{i'})$ for $1\leq i<i'\leq 3$.
By the symmetry, we may assume that $V(T)\cap N_{G}(w_{1})=V(T)\cap N_{G}(w_{2})=\{v_{1},v_{2}\}$ and $V(T)\setminus N_{G}(w_{3})=\{v_{j}\}$ for $j\in \{2,3\}$.
Now we consider the mapping $c':V(G)\rightarrow [3]$ with
$$
c'(a)=
\begin{cases}
1 & (a\in \{u,v_{1}\})\\
2 & (a=v_{2})\\
3 & (a\in \{v_{3},w_{1},w_{2}\})\\
j & (a=w_{3}).
\end{cases}
$$
Then we can easily verify that $c'$ is a proper $3$-coloring of $G$, which is a contradiction.
Thus $V(T)\cap N_{G}(w_{i})\neq V(T)\cap N_{G}(w_{i'})$ for all $1\leq i<i'\leq 3$.
This implies that $G=H_{0}$.
\qed

\begin{lem}%%%%%%%%%%%%%%%%%%%%%%%%%%%%%%%%%%%%%%%%%%%%%%%%%%%%%%%%%%%%%%%%%%%%%%%%%%%%%%%%%%%%%%%%%%%%%%%%%%%%%%%%%%%%%
\label{lem-small-3}
Let $k\geq 5$ be an integer, and let $n$ be a positive integer.
Then all graphs $G$ and $H$ with $|V(G)|\leq k+3$, $|V(H)|\leq n$ and $\min\{\chi (G),\chi (H)\}=k$ satisfy $\chi (G\times H)=k$ if and only if all $K_{k-1}$-free $k$-critical graphs $H$ with $k+4\leq |V(H)|\leq n$ satisfy $\chi ((K_{k-4}+H_{0})\times H)=k$.
\end{lem}
%%%%%%%%%%%%%%%%%%%%%%%%%%%%%%%%%%%%%%%%%%%%%%%%%%%%%%%%%%%%%%%%%%%%%%%%%%%%%%%%%%%%%%%%%%%%%%%%%%%%%%%%%%%%%%%%%%%%%%%%
\proof
The ``only if'' part is trivial. Thus we show the ``if'' part.
We suppose that
\begin{align}
\mbox{all $K_{k-1}$-free $k$-critical graphs $H$ with $k+4 \leq |V(H)| \leq n$ satisfy $\chi ((K_{k-4}+H_{0})\times H)=k$.}\label{cond-lem-small-3-1}
\end{align}
Let $G'$ and $H'$ be $k$-critical subgraphs of $G$ and $H$, respectively.
Since $G'\times H'$ is a subgraph of $G\times H$, we have $\chi (G\times H)\geq \chi (G'\times H')$.
Considering (\ref{eq-intro-1}), it suffices to show that $\chi (G'\times H')=k$.

We first assume that $\min\{|V(G')|,|V(H')|\}\leq k+2$.
Without loss of generality, we may assume that $|V(G')|\leq k+2$.
Then by Lemma~\ref{lem-small-1}, $G'$ is either $K_{k}$ or $K_{k-3}+C_{5}$.
In particular, each vertex of $G'$ belongs to a $(k-1)$-clique of $G'$.
Hence by Theorem~\ref{Thm-prel-B}, $\chi (G'\times H')=k$, as desired.
Thus we may assume that $|V(G')|=k+3$ and $|V(H')|\geq k+3$.

By Lemma~\ref{lem-small-1}, $G'=K_{k-4}+A$ for some $A\in \AA_{4}$.
Suppose that $A\neq H_{0}$.
Then by Lemma~\ref{lem-small-2}, each vertex of $A$ belongs to a triangle.
Since $G'=K_{k-4}+A$, this implies that each vertex of $G'$ belongs to a $(k-1)$-clique of $G'$.
This together with Theorem~\ref{Thm-prel-B} implies that $\chi (G'\times H')=k$.
Thus we may assume that $A=H_{0}$.

If $H'$ contains $(k-1)$-clique, then both $G'$ and $H'$ contain $(k-1)$-cliques, and hence by Theorem~\ref{Thm-prel-C}, $\chi (G'\times H')=k$, as desired.
Thus we may assume that $H'$ is $K_{k-1}$-free.
If $|V(H')|=k+3$, then by similar argument in the previous paragraph, we have $H'=K_{k-4}+H_{0}$, which contradicts the $K_{k-1}$-freeness of $H'$.
Thus $k+4 \leq |V(H')|\leq |V(H)| \leq n$. 
Then by (\ref{cond-lem-small-3-1}), $\chi (G'\times H')=k$.
\qed

Now we prove Theorem~\ref{thm-small-main}.

\medbreak\noindent\textit{Proof of Theorem~\ref{thm-small-main}.}\quad
Applying Lemma~\ref{lem-small-3} with $n=k+3$, we obtain that if one of the assumptions (i) and (ii) of the theorem holds, then $\chi (G\times H)=k$.
Thus we may assume that $k=5$ and $(|V(G)|,|V(H)|)\in \{(8,8),(8,9),(8,10),(9,8),(10,8)\}$.
We may assume that $|V(G)|=8$.
Then by Lemma~\ref{lem-small-3} with $n=10$, it suffices to show that all $K_{4}$-free $5$-critical graphs $H'$ with $9\leq |V(H)|\leq 10$ satisfy $\chi ((K_{1}+H_{0})\times H)=5$.
However, it follows from Theorem~\ref{Thm-smallcolor} that every $K_{4}$-free $5$-critical graph has at least $11$ vertices, and so there is no target graph.

This completes the proof of Theorem~\ref{thm-small-main}.
\qed

To prove Theorem~\ref{thm-small-main2}, we prepare the following lemma.

\begin{lem}%%%%%%%%%%%%%%%%%%%%%%%%%%%%%%%%%%%%%%%%%%%%%%%%%%%%%%%%%%%%%%%%%%%%%%%%%%%%%%%%%%%%%%%%%%%%%%%%%%%%%%%%%%%%%
\label{lem-small-k>=6-1}
For an integer $k\geq 6$, a $K_{k-1}$-free graph $G$ of order $k+4$ is $k$-critical if and only if $G=K_{k-6}+C_{5}+C_{5}$.
\end{lem}
%%%%%%%%%%%%%%%%%%%%%%%%%%%%%%%%%%%%%%%%%%%%%%%%%%%%%%%%%%%%%%%%%%%%%%%%%%%%%%%%%%%%%%%%%%%%%%%%%%%%%%%%%%%%%%%%%%%%%%%%
\proof
The ``if'' part is trivial.
Thus we show the ``only if'' part by induction on $k$.
Note that $|V(G)|=k+4$ and $k+4<2k-1$ (i.e. $k\geq 6$) holds by our assumption.
Hence, we have $|V(G)|=k+4<2k-1$.
Then $G$ is decomposable by Theorem~\ref{Thm-prel-A}, and hence $G=G_{1}+G_{2}$ for a $k_{1}$-critical graph $G_{1}$ and a $k_{2}$-critical graph $G_{2}$ with $k_{1}+k_{2}=k$ and $k_{1}\geq k_{2}$.
Choose $G_{1}$ and $G_{2}$ so that $k_{2}$ is as small as possible.

For the moment, we suppose that $G_{1}$ and $G_{2}$ are non-complete.
Then by Lemma~\ref{lem-small-00}, $|V(G_{i})|\geq k_{i}+2$ for each $i\in \{1,2\}$.
Since
$$
k+4=|V(G)|=|V(G_{1})|+|V(G_{2})|\geq (k_{1}+2)+(k_{2}+2)=k+4,
$$
we have $|V(G_{i})|=k_{i}+2$.
It follows from Lemma~\ref{lem-small-1} that $G_1=K_{k_1-3}+C_5$ and $G_2=K_{k_2-3}+C_5$. 
Hence, $G=G_1+G_2=K_{k_1+k_2-6}+C_5+C_5$, as desired. 

If $k=6$, then $G_{1}$ is a $K_{4}$-free $5$-critical graph of order $|V(G)\setminus V(G_{2})|~(=9)$, which contradicts Theorem~\ref{Thm-smallcolor}.
Thus $k\geq 7$ (and the first step of the induction is completed).
Since $G_{1}$ is a $K_{k-2}$-free $(k-1)$-critical graph of order $k+3$, we have $G_{1}=K_{(k-1)-6}+C_{5}+C_{5}$ by the induction hypothesis.
Consequently, $G=G_{1}+K_{1}=K_{k-6}+C_{5}+C_{5}$, as desired.
\qed

Combining Lemmas~\ref{lem-small-3} and \ref{lem-small-k>=6-1}, we obtain Theorem~\ref{thm-small-main2}.

%%%%%%%%%%%%%%%%%%%%%%%%%%%%%%%%%%%%%%%%%%%%%%%%%%%%%%%%%%%%%%%%%%%%%%%%%%%%%%%%%%%%%%%%%%%%%%%%%%%%%%%%%%%%%%%%%%%%%%%%
%%%%%%%%%%%%%%%%%%%%%%%%%%%%%%%%%%%%%%%%%%%%%%%%%%%%%%%%%%%%%%%%%%%%%%%%%%%%%%%%%%%%%%%%%%%%%%%%%%%%%%%%%%%%%%%%%%%%%%%%
%%%%%%%%%%%%%%%%%%%%%%%%%%%%%%%%%%%%%%%%%%%%%%%%%%%%%%%%%%%%%%%%%%%%%%%%%%%%%%%%%%%%%%%%%%%%%%%%%%%%%%%%%%%%%%%%%%%%%%%%
\section{Algebraic reduction of Conjecture~\ref{conj1}}\label{sec-main}
%%%%%%%%%%%%%%%%%%%%%%%%%%%%%%%%%%%%%%%%%%%%%%%%%%%%%%%%%%%%%%%%%%%%%%%%%%%%%%%%%%%%%%%%%%%%%%%%%%%%%%%%%%%%%%%%%%%%%%%%
%%%%%%%%%%%%%%%%%%%%%%%%%%%%%%%%%%%%%%%%%%%%%%%%%%%%%%%%%%%%%%%%%%%%%%%%%%%%%%%%%%%%%%%%%%%%%%%%%%%%%%%%%%%%%%%%%%%%%%%%
%%%%%%%%%%%%%%%%%%%%%%%%%%%%%%%%%%%%%%%%%%%%%%%%%%%%%%%%%%%%%%%%%%%%%%%%%%%%%%%%%%%%%%%%%%%%%%%%%%%%%%%%%%%%%%%%%%%%%%%%

%%%%%%%%%%%%%%%%%%%%%%%%%%%%%%%%%%%%%%%%%%%%%%%%%%%%%%%%%%%%%%%%%%%%%%%%%%%%%%%%%%%%%%%%%%%%%%%%%%%%%%%%%%%%%%%%%%%%%%%%
%%%%%%%%%%%%%%%%%%%%%%%%%%%%%%%%%%%%%%%%%%%%%%%%%%%%%%%%%%%%%%%%%%%%%%%%%%%%%%%%%%%%%%%%%%%%%%%%%%%%%%%%%%%%%%%%%%%%%%%%
%%%%%%%%%%%%%%%%%%%%%%%%%%%%%%%%%%%%%%%%%%%%%%%%%%%%%%%%%%%%%%%%%%%%%%%%%%%%%%%%%%%%%%%%%%%%%%%%%%%%%%%%%%%%%%%%%%%%%%%%
\subsection{Equivalence conjecture for Conjecture~\ref{conj1} via the criticality}\label{sec-main1}
%%%%%%%%%%%%%%%%%%%%%%%%%%%%%%%%%%%%%%%%%%%%%%%%%%%%%%%%%%%%%%%%%%%%%%%%%%%%%%%%%%%%%%%%%%%%%%%%%%%%%%%%%%%%%%%%%%%%%%%%
%%%%%%%%%%%%%%%%%%%%%%%%%%%%%%%%%%%%%%%%%%%%%%%%%%%%%%%%%%%%%%%%%%%%%%%%%%%%%%%%%%%%%%%%%%%%%%%%%%%%%%%%%%%%%%%%%%%%%%%%
%%%%%%%%%%%%%%%%%%%%%%%%%%%%%%%%%%%%%%%%%%%%%%%%%%%%%%%%%%%%%%%%%%%%%%%%%%%%%%%%%%%%%%%%%%%%%%%%%%%%%%%%%%%%%%%%%%%%%%%%

In this subsection, we focus on the following conditions for given graphs $G$ and $H$:
\vspace{-.5\baselineskip}
\begin{enumerate}
\setlength{\parskip}{0cm}
\setlength{\itemsep}{0cm}
\item[{\bf (X1)}]
$\chi (G\times H)\leq k-1$;
\item[{\bf (W1)}]
$\min\{\chi (G),\chi (H)\}\leq k-1$;
\item[{\bf (V1)}]
$\delta (G)\geq 1$ and there exists a proper $(k-1)$-coloring $c$ of $G-uv$ with $c(u)=c(v)=1$ for all $uv\in E(G)$;
\item[{\bf (V2)}]
$\delta (H)\geq 1$ and there exists a proper $(k-1)$-coloring $c'$ of $H-u'v'$ with $c'(u')=c'(v')=1$ for all $u'v'\in E(H)$;
\item[{\bf (V3)}]
there exists a vertex of $G$ belonging to no $(k-1)$-clique of $G$;
\item[{\bf (V4)}]
there exists a vertex of $H$ belonging to no $(k-1)$-clique of $H$;
\item[{\bf (V5)}]
$\max\{\omega (G),\omega (H)\}\leq k-1$ and $\min\{\omega (G),\omega (H)\}\leq k-2$;
\item[{\bf (V6)}]
$\delta (G)\geq k-1$; and
\item[{\bf (V7)}]
$\delta (H)\geq k-1$.
\end{enumerate}
\vspace{-.5\baselineskip}

Note that the conditions (V1)--(V7) derive from the definitions or the previous results in Section~\ref{sec-prel} as follows: 
\vspace{-.5\baselineskip}
\begin{enumerate}[$\bullet $]
\setlength{\parskip}{0cm}
\setlength{\itemsep}{0cm}
\item (V1) and (V2) derive from the definition of the criticallity. 
\item If (V3) or (V4) is not satisfied for graphs $G$ and $H$ with $\chi(G)=\chi(H)=k$, then Conjecture~\ref{conj1} is automatically true by Theorem~\ref{Thm-prel-B}. 
\item Assume $G$ and $H$ are $k$-critical graphs.
On (V5), if $\max\{\omega (G),\omega (H)\}\geq k$, say, $\omega(G) \geq k$, 
then we see that $G$ should be a complete graph of order $k$ by the criticality of $G$, so Conjecture~\ref{conj1} is automatically true by Theorem~\ref{Thm-prel-B}. 
If $\min\{\omega (G),\omega (H)\} \geq k-1$, then Conjecture~\ref{conj1} is automatically true by Theorem~\ref{Thm-prel-C}. 
\item (V6) and (V7) always hold for $k$-critical graphs $G$ and $H$, respectively, by Theorem~\ref{prop-verc-mindeg}. 
\end{enumerate}

For an integer $n\geq 1$, let $\GG_{n}$ be the set of graphs of order at most $n$.
We define two sets as follows:
\begin{align*}
W_{k,n,n'} &= \{(G,H) \in \GG_{n}\times \GG_{n'} : (G,H)\text{ satisfies } \text{(X1) and (W1)}\}; \mbox{ and}\\
V_{k,n,n'} &= \{(G,H) \in \GG_{n}\times \GG_{n'} : (G,H)\text{ satisfies } \text{(X1) and (V1)--(V7)}\}.
\end{align*}

The following is the key proposition for our argument.

\begin{prop}%%%%%%%%%%%%%%%%%%%%%%%%%%%%%%%%%%%%%%%%%%%%%%%%%%%%%%%%%%%%%%%%%%%%%%%%%%%%%%%%%%%%%%%%%%%%%%%%%%%%%%%%%%%%
\label{prop-main1-1}
For integers $k\geq 2$, $n\geq 1$ and $n'\geq 1$, the following are equivalent:
\vspace*{-2.5mm}
\begin{enumerate}[{\bf (H1)}]
\setlength{\parskip}{0cm}
\setlength{\itemsep}{0cm}
\item
$V_{k,n,n'}\subseteq W_{k,n,n'}$;
\item
if $G\in \GG_{n}$ and $H\in \GG_{n'}$ are $k$-critical, then $\chi (G\times H) = k$; and
\item
if $G\in \GG_{n}$ and $H\in \GG_{n'}$ satisfy $\min\{\chi (G),\chi (H)\}=k$, then $\chi (G\times H)=k$.
\end{enumerate}
\vspace*{-2.5mm}
\end{prop}
%%%%%%%%%%%%%%%%%%%%%%%%%%%%%%%%%%%%%%%%%%%%%%%%%%%%%%%%%%%%%%%%%%%%%%%%%%%%%%%%%%%%%%%%%%%%%%%%%%%%%%%%%%%%%%%%%%%%%%%%
\proof
%``(H1)~$\Rightarrow $~(H2)'' is clear because $V'_{k,n,n'}\subseteq V_{k,n,n'}$ trivially holds by definition.
We first prove ``(H1)~$\Rightarrow $~(H2)''.
Suppose that (H1) holds and there exist $k$-critical graphs $G\in \GG_{n}$ and $H\in \GG_{n'}$ such that $\chi (G\times H)\leq k-1$ (i.e., (X1) holds).
By the $k$-criticality of $G$ and $H$, (V1) and (V2) clearly hold.
If each vertex of $G$ belongs to a $(k-1)$-clique of $G$, then by Theorem~\ref{Thm-prel-B}, $\chi (G\times H)=k$, a contradiction. 
Thus both (V3) and (V4) hold. If one of $G$ and $H$ contains a $k$-clique, then by its $k$-criticality, it is a complete graph of order $k$, which contradicts (V3) or (V4).
Thus $\max\{\omega (G),\omega (H)\}\leq k-1$.
If $\min\{\omega (G),\omega (H)\}\geq k-1$, then by Theorem~\ref{Thm-prel-C}, $\chi (G\times H)=k$, a contradiction.
Therefore (V5) holds.
Furthermore, it follows from Theorem~\ref{prop-verc-mindeg} that (V6) and (V7) hold. Consequently, we have $(G,H)\in V_{k,n,n'}$.
By our assumption, we have $(G,H)\in W_{k,n,n'}$. In particular, $\min\{\chi (G),\chi (H)\}\leq k-1$, which contradicts the assumption that $G$ and $H$ are $k$-critical. 

We next prove ``(H2)~$\Rightarrow $~(H1)''.
Suppose that (H2) holds.
Let $(G,H)\in V_{k,n,n'}$.
Then $G$ and $H$ satisfy (X1) and (V1)--(V7).
If $\chi (G)\geq k$ and $\chi (H)\geq k$, then the conditions (V1) and (V2) force both $G$ and $H$ to be $k$-critical, and hence $\chi (G\times H)=k$ by (H2), which contradicts (X1).
Thus $\chi (G)\leq k-1$ or $\chi (H)\leq k-1$.
In particular, $G$ and $H$ satisfy (W1), and so $(G,H)\in W_{k,n,n'}$. 

Finally, we prove ``(H2)~$\Leftrightarrow $~(H3)''.
Since ``(H3)~$\Rightarrow $~(H2)'' trivially holds, it suffices to show that ``(H2)~$\Rightarrow $~(H3)'' holds.
Suppose that (H2) holds.
Let $G\in \GG_{n}$ and $H\in \GG_{n'}$ be graphs with $\min\{\chi (G),\chi (H)\}=k$.
Then $G$ contains a $k$-critical subgraph $G'$ and $H$ contains a $k$-critical subgraph $H'$.
By (H2), we have $\chi (G'\times H')=k$, and hence $\chi (G\times H)\geq \chi (G'\times H')=k$ because $G'\times H'$ is a subgraph of $G\times H$.
This together with (\ref{eq-intro-1}) implies that $\chi (G\times H)=k$.

This completes the proof of the proposition.
\qed

By Theorem~\ref{thm-small-main} and Proposition~\ref{prop-main1-1}, we can translate Conjecture~\ref{conj1} into an inclusion relation problem concerning $W_{k,n,n'}$ and $V_{k,n,n'}$ as follows.

\begin{cor}%%%%%%%%%%%%%%%%%%%%%%%%%%%%%%%%%%%%%%%%%%%%%%%%%%%%%%%%%%%%%%%%%%%%%%%%%%%%%%%%%%%%%%%%%%%%%%%%%%%%%%%%%%%%%
\label{cor2-1-2}
Let $k\geq 3$, $n\geq 1$ and $n'\geq 1$ be integers. Then the following are equivalent:
\begin{enumerate}[{\upshape(i)}]
\item
Conjecture~\ref{conj1} is true for the case where $\min\{\chi (G),\chi (H)\}=k$;
\item
$V_{k,n,n'}\subseteq W_{k,n,n'}$ for any integers $n\geq k+3$ and $n'\geq k+3$ with $(n,n')\neq (k+3,k+3)$. 
\end{enumerate}
\end{cor}
%%%%%%%%%%%%%%%%%%%%%%%%%%%%%%%%%%%%%%%%%%%%%%%%%%%%%%%%%%%%%%%%%%%%%%%%%%%%%%%%%%%%%%%%%%%%%%%%%%%%%%%%%%%%%%%%%%%%%%%%

%%%%%%%%%%%%%%%%%%%%%%%%%%%%%%%%%%%%%%%%%%%%%%%%%%%%%%%%%%%%%%%%%%%%%%%%%%%%%%%%%%%%%%%%%%%%%%%%%%%%%%%%%%%%%%%%%%%%%%%%
%%%%%%%%%%%%%%%%%%%%%%%%%%%%%%%%%%%%%%%%%%%%%%%%%%%%%%%%%%%%%%%%%%%%%%%%%%%%%%%%%%%%%%%%%%%%%%%%%%%%%%%%%%%%%%%%%%%%%%%%
%%%%%%%%%%%%%%%%%%%%%%%%%%%%%%%%%%%%%%%%%%%%%%%%%%%%%%%%%%%%%%%%%%%%%%%%%%%%%%%%%%%%%%%%%%%%%%%%%%%%%%%%%%%%%%%%%%%%%%%%
\subsection{Equivalence conjecture for Conjecture~\ref{conj1} via ideals of polynomial rings}\label{sec-main2}
%%%%%%%%%%%%%%%%%%%%%%%%%%%%%%%%%%%%%%%%%%%%%%%%%%%%%%%%%%%%%%%%%%%%%%%%%%%%%%%%%%%%%%%%%%%%%%%%%%%%%%%%%%%%%%%%%%%%%%%%
%%%%%%%%%%%%%%%%%%%%%%%%%%%%%%%%%%%%%%%%%%%%%%%%%%%%%%%%%%%%%%%%%%%%%%%%%%%%%%%%%%%%%%%%%%%%%%%%%%%%%%%%%%%%%%%%%%%%%%%%
%%%%%%%%%%%%%%%%%%%%%%%%%%%%%%%%%%%%%%%%%%%%%%%%%%%%%%%%%%%%%%%%%%%%%%%%%%%%%%%%%%%%%%%%%%%%%%%%%%%%%%%%%%%%%%%%%%%%%%%%

Throughout this section, we fix integers $k\geq 3$, $n\geq 1$ and $n'\geq 1$.
We start with an easy algebraic proposition.

\begin{prop}%%%%%%%%%%%%%%%%%%%%%%%%%%%%%%%%%%%%%%%%%%%%%%%%%%%%%%%%%%%%%%%%%%%%%%%%%%%%%%%%%%%%%%%%%%%%%%%%%%%%%%%%%%%%
\label{prop-main2-1}
Let $e$, $x_{1}$ and $x_{2}$ be three variables satisfying $e(e-1)=0$ and $x_{1}^{k}-1=x_{2}^{k}-1=0$. 
Then $e(x_{1}^{k-1}+x_{1}^{k-2}x_{2}+\cdots +x_{2}^{k-1})=0$ if and only if $e=0$ or $x_{1}\neq x_{2}$.
\end{prop}
%%%%%%%%%%%%%%%%%%%%%%%%%%%%%%%%%%%%%%%%%%%%%%%%%%%%%%%%%%%%%%%%%%%%%%%%%%%%%%%%%%%%%%%%%%%%%%%%%%%%%%%%%%%%%%%%%%%%%%%%
\proof
Note that $x_{i}\neq 0$ for each $i\in \{1,2\}$.
If $x_{1}\neq x_{2}$, then $x_{1}^{k-1}+x_{1}^{k-2}x_{2}+\cdots +x_{2}^{k-1}=0$ because
$$
0=(x_{1}^{k}-1)-(x_{2}^{k}-1)=(x_{1}-x_2)(x_{1}^{k-1}+x_{1}^{k-2}x_{2}+\cdots +x_{2}^{k-1}).
$$
Conversely, if $x_{1}=x_{2}$, then $x_{1}^{k-1}+x_{1}^{k-2}x_{2}+\cdots +x_{2}^{k-1}=kx_{1}^{k-1}\neq 0$.
Hence $x_{1}\neq x_{2}$ if and only if $x_{1}^{k-1}+x_{1}^{k-2}x_{2}+\cdots +x_{2}^{k-1}=0$, which proves the proposition.
\qed

We prepare the variables
$$
x_{1},\ldots ,x_{n},e_{ij}~(1\leq i<j\leq n),~~~y_{1},\ldots ,y_{n'},f_{i'j'}~(1\leq i'<j'\leq n'),~~~z_{ii'}~(1\leq i\leq n,~1\leq i'\leq n').
$$
Then, considering Proposition~\ref{prop-main2-1}, we obtain the following: The solutions of system of equations 
$$
\begin{cases}
e_{ij}(e_{ij}-1)=0~~~(1\leq i<j\leq n)\\
x_{i}^{k}-1=0~~~(1\leq i\leq n)\\
e_{ij}(x_{i}^{k-1}+x_{i}^{k-2}x_{j}+\cdots +x_{j}^{k-1})=0~~~(1\leq i<j\leq n)
\end{cases}
$$
one-to-one correspond to the pairs of a labeled graph on $[n]$ and its proper $k$-coloring; the solutions of system of equations 
$$
\begin{cases}
f_{i'j'}(f_{i'j'}-1)=0~~~(1\leq i'<j'\leq n')\\
y_{i'}^{k}-1=0~~~(1\leq i'\leq n')\\
f_{i'j'}(y_{i'}^{k-1}+y_{i'}^{k-2}y_{j'}+\cdots +y_{j'}^{k-1})=0~~~(1\leq i'<j'\leq n')
\end{cases}
$$
one-to-one correspond to the pairs of a labeled graph on $[n']$ and its proper $k$-coloring; and the solutions of system of equations 
$$
\begin{cases}
e_{ij}(e_{ij}-1)=0~~~(1\leq i<j\leq n)\\
f_{i'j'}(f_{i'j'}-1)=0~~~(1\leq i'<j'\leq n')\\
z_{ii'}^{k}-1=0~~~(1\leq i\leq n,~1\leq i'\leq n')\\
e_{ij}f_{i'j'}(z_{ii'}^{k-1}+z_{ii'}^{k-2}z_{jj'}+\cdots +z_{jj'}^{k-1})=0~~~(1\leq i<j\leq n,~1\leq i'<j'\leq n')
\end{cases}
$$
one-to-one correspond to the pairs of the tensor product of labeled graphs on $[n]$ and $[n']$ and its proper $k$-coloring.
We explain an outline of, for example, the first fact.
Consider a graph on $[n]$.
We regard a solution of $e_{ij}(e_{ij}-1)=0~(1\leq i<j\leq n)$ as its adjacency matrix, and a solution of $x_{i}^{k}-1=0$ as a color assigned to the vertex $i$ (here $x_{i}$ can take exactly $k$ solutions because $x_{i}$ is a $k$-th root of unity).
Then $e_{ij}(x_{i}^{k-1}+x_{i}^{k-2}x_{j}+\cdots +x_{j}^{k-1})=0~(1\leq i<j\leq n)$ implies by Proposition~\ref{prop-main2-1} that 
if two vertices $i$ and $j$ are adjacent, then the color assigned to $i$ differs from the color assigned to $j$.

Based on the above facts, we associate solutions of some systems of equations with the members in $W_{k,n,n'}$ and $V_{k,n,n'}$ appearing in Subsection~\ref{sec-main1}.

\subsection*{Description of $W_{k,n,n'}$}
All the ideals below (i.e., the ideals $E_{n,n'}$, $X_{k,n,n'}$, $Z_{k,n,n'}$, $I_{k,n}$, $I_{k,n'}'$ and $J_{k,n,n'}$) are regarded 
as the ones of the polynomial ring $$\mathbb{C}[e_{ij},f_{i'j'},x_s,y_{s'},z_{ss'} : 1 \leq i<j \leq n, 1 \leq i'<j' \leq n', 1 \leq s \leq n, 1 \leq s' \leq n']$$ of $(\binom{n}{2}+\binom{n'}{2}+n+n'+nn')$ variables. 

We define several ideals as follows: 
\begin{align}
E_{n,n'}&=(e_{ij}(e_{ij}-1):1\leq i<j\leq n)+(f_{i'j'}(f_{i'j'}-1):1\leq i'<j'\leq n'), \tag{$E(G)$ and $E(H)$} \\
X_{k,n,n'}&=(x_{i}^{k-1}-1:1\leq i\leq n)+(y_{i'}^{k-1}-1:1\leq i'\leq n'), \tag{$(k-1)$-colorings of $G$ and $H$} \\
Z_{k,n,n'}&=(z_{ii'}^{k-1}-1:1\leq i\leq n,1\leq i'\leq n'), \tag{$(k-1)$-coloring of $G\times H$} \\
I_{k,n} &= (e_{ij}(x_{i}^{k-2}+x_{i}^{k-3}x_{j}+\cdots +x_{j}^{k-2}):1\leq i<j\leq n), \tag{W1($G$)} \\
I'_{k,n'} &= (f_{i'j'}(y_{i'}^{k-2}+y_{i'}^{k-3}y_{j'}+\cdots +y_{j'}^{k-2}):1\leq i'<j'\leq n'). \tag{W1($H$)} 
\end{align}
Note that the solutions of $E_{n,n'}+X_{k,n,n'}+I_{k,n} \cdot I'_{k,n'}$ 
one-to-one correspond to the pairs of graphs $(G,H) \in \GG_n \times \GG_{n'}$ satisfying (W1) with their proper $(k-1)$-colorings. 
Furthermore, let
\begin{align}\tag{X1}
J_{k,n,n'} &= (e_{ij}f_{i'j'}(z_{ii'}^{k-2}+z_{ii'}^{k-3}z_{jj'}+\cdots +z_{jj'}^{k-2}):1\leq i<j\leq n,1\leq i'<j'\leq n'). 
\end{align}
Let
$$ \JJ_{k,n,n'}=E_{n,n'}+X_{k,n,n'}+Z_{k,n,n'}+I_{k,n} \cdot I'_{k,n'}+J_{k,n,n'} $$
be the ideal, and set 
\begin{align}\tag{$W_{k,n,n'}$}
\tilde{\JJ}_{k,n,n'}=\JJ_{k,n,n'}\cap \mathbb{C}[e_{ij},f_{i'j'}].
\end{align}
Then we can verify that the solutions of $\tilde{\JJ}_{k,n,n'}$ one-to-one correspond to the members of $W_{k,n,n'}$.

\subsection*{Description of $V_{k,n,n'}$} 
All the ideals below (i.e., the ideals appearing in $\II_{k,n,n'}$) are regarded as the ideals of the polynomial ring 
\begin{align*}
\mathbb{C}[e_{ij},f_{i'j'},x_{pq\ell},y_{p'q'\ell'},z_{ss'} : &1 \leq i<j \leq n, 1 \leq i'<j' \leq n', 1 \leq p < q \leq n, 1 \leq p'<q' \leq n', \\
&1 \leq \ell \leq n, 1 \leq \ell' \leq n', 1 \leq s \leq n, 1 \leq s' \leq n']
\end{align*} 
of $(\binom{n}{2}(n+1)+\binom{n'}{2}(n'+1)+nn')$ variables.

We define several ideals as follows: 
\begin{align}\tag{V1}
\begin{split}
P_{k,n} &= ((\prod _{1\leq j<i}(e_{ji}-1))\cdot (\prod _{i<j\leq n}(e_{ij}-1)):1\leq i\leq n)\\
&+ (x_{pqi}^{k-1}-1:1\leq p<q\leq n,~1\leq i\leq n)+(x_{pqi}-1:1\leq p<q\leq n,~i\in \{p,q\})\\
&+ (e_{pq}e_{ij}(x_{pqi}^{k-2}+x_{pqi}^{k-3}x_{pqj}+\cdots +x_{pqj}^{k-2}):1\leq p<q\leq n,~1\leq i<j\leq n,~(p,q)\neq (i,j)), 
\end{split}
\end{align}
\begin{align}\tag{V2}
\begin{split}
P'_{k,n'} &= ((\prod _{1\leq j'<i'}(f_{j'i'}-1))\cdot (\prod _{i'<j'\leq n'}(f_{i'j'}-1)):1\leq i'\leq n')\\
&+ (y_{p'q'i'}^{k-1}-1:1\leq p'<q'\leq n',~1\leq i'\leq n')+(y_{p'q'i'}-1:1\leq p'<q'\leq n',~i'\in \{p',q'\})\\
&+ (f_{p'q'}f_{i'j'}(y_{p'q'i'}^{k-2}+y_{p'q'i'}^{k-3}y_{p'q'j'}+\cdots +y_{p'q'j'}^{k-2}):1\leq p'<q'\leq n',~1\leq i'<j'\leq n',~(p',q')\neq (i',j')), 
\end{split}
\end{align}
\begin{align}
Q_{k,n} &= ((\prod_{i\in X}e_{1i})\cdot (\prod_{\substack{i,j\in X\\ i<j}}e_{ij}):X\subseteq [n]\setminus \{1\}\mbox{ with }|X|=k-2), \tag{V3} \\
Q'_{k,n'} &= ((\prod_{i'\in X'}f_{1i'})\cdot (\prod_{\substack{i',j'\in X'\\ i'<j'}}f_{i'j'}):X'\subseteq [n']\setminus \{1\}\mbox{ with }|X'|=k-2), \tag{V4} \\
R_{k,n} &= (\prod_{\substack{i,j\in X\\ i<j}}e_{ij}:X\subseteq [n]\mbox{ with }|X|=k), \tag{$\omega(G) \leq k-1$} \\
R'_{k,n'} &= (\prod_{\substack{i',j'\in X'\\ i'<j'}}f_{i'j'}:X'\subseteq [n']\mbox{ with }|X'|=k). \tag{$\omega(H) \leq k-1$} 
\end{align}
Note that the condition that $\omega(G) \leq k-1$ {\em and} $\omega(H) \leq k-1$ hold is equivalent to $\max\{\omega(G),\omega(H)\} \leq k-1$, 
while the condition that $\omega(G) \leq k-2$ {\em or} $\omega(H) \leq k-2$ holds is equivalent to $\min\{\omega(G),\omega(H)\} \leq k-2$. 
Let \begin{align}
S_{k,n}&=((\prod_{\substack{i\in X\\ i<\ell}}(e_{i\ell}-1))\cdot (\prod_{\substack{i\in X\\ \ell<i}}(e_{\ell i}-1)):\ell\in [n],~X\subseteq [n]\setminus \{\ell\}\mbox{ with }|X|=n-k+1), \tag{V6} \\
S'_{k,n'}&=((\prod_{\substack{i'\in X'\\ i'<\ell'}}(f_{i'\ell'}-1))\cdot (\prod_{\substack{i'\in X'\\ \ell'<i'}}(f_{\ell'i'}-1)):\ell'\in [n'],~X'\subseteq [n']\setminus \{\ell'\}\mbox{ with }|X'|=n'-k+1). \tag{V7}
\end{align}
Furthermore, let 
\begin{align*}
\II_{k,n,n'}=E_{n,n'}&+Z_{k,n,n'}+J_{k,n,n'}+P_{k,n}+P'_{k,n'}+Q_{k,n}+Q'_{k,n'}\\
&+\underbrace{R_{k,n}+R'_{k,n'}+R_{k-1,n}\cdot R'_{k-1,n'}}_\text{(V5)}+S_{k,n}+S'_{k,n'} 
\end{align*}
and set
\begin{align}\tag{$V_{k,n,n'}$}
\tilde{\II}_{k,n,n'}=\II_{k,n,n'}\cap \mathbb{C}[e_{ij},f_{i'j'}].
\end{align}
Then we can verify that the solutions of $\tilde{\II}_{k,n,n'}$ one-to-one correspond to the members of $V_{k,n,n'}$.

%%%%%%%%%%%%%%%%%%%%%%%%%%%%%%%%%%%%%%%%%%%%%%%%%%%%%%%%%%%%%%%%%%%%%%%%%%%%%%%%%%%%%%%%%%%%%%%%%%%%%%%%%%%%%%%%%%%%%%%%%%%%%%%%%%%%%%%%%%%%%%%%%
%%%%%%%%%%%%%%%%%%%%%%%%%%%%%%%%%%%%%%%%%%%%%%%%%%%%%%%%%%%%%%%%%%%%%%%%%%%%%%%%%%%%%%%%%%%%%%%%%%%%%%%%%%%%%%%%%%%%%%%%%%%%%%%%%%%%%%%%%%%%%%%%%
%%%%%%%%%%%%%%%%%%%%%%%%%%%%%%%%%%%%%%%%%%%%%%%%%%%%%%%%%%%%%%%%%%%%%%%%%%%%%%%%%%%%%%%%%%%%%%%%%%%%%%%%%%%%%%%%%%%%%%%%%%%%%%%%%%%%%%%%%%%%%%%%%
%%%%%%%%%%%%%%%%%%%%%%%%%%%%%%%%%%%%%%%%%%%%%%%%%%%%%%%%%%%%%%%%%%%%%%%%%%%%%%%%%%%%%%%%%%%%%%%%%%%%%%%%%%%%%%%%%%%%%%%%%%%%%%%%%%%%%%%%%%%%%%%%%

\bigskip

Consequently, it follows from Proposition~\ref{prop-main1-1} that the following theorem holds. 

\begin{thm}%%%%%%%%%%%%%%%%%%%%%%%%%%%%%%%%%%%%%%%%%%%%%%%%%%%%%%%%%%%%%%%%%%%%%%%%%%%%%%%%%%%%%%%%%%%%%%%%%%%%%%%%%%%%%
\label{thm3-1}
Let $k\geq 3$, $n\geq 1$ and $n'\geq 1$ be integers.
Then the following are equivalent:
\begin{enumerate}[{\upshape(i)}]
\item
Conjecture~\ref{conj1} is true for the case where $G\in \GG_{n}$, $H\in \GG_{n'}$ and $\min\{\chi (G),\chi (H)\}=k$; 
\item
$\tilde{\JJ}_{k,n,n'}\subseteq \tilde{\II}_{k,n,n'}$.
%\item
%$\tilde{\JJ}_{k,n,n'}\subseteq \tilde{\LL}_{k,n,n'}$.
\end{enumerate}
\end{thm}
%%%%%%%%%%%%%%%%%%%%%%%%%%%%%%%%%%%%%%%%%%%%%%%%%%%%%%%%%%%%%%%%%%%%%%%%%%%%%%%%%%%%%%%%%%%%%%%%%%%%%%%%%%%%%%%%%%%%%%%%

\begin{rem}{\em 
For the computations of $\tilde{\JJ}_{k,n,n'}$ and $\tilde{\II}_{k,n,n'}$, we have to \textit{eliminate} the variables. 
For example, $\tilde{\JJ}_{k,n,n'}$ is defined by $\JJ_{k,n,n'} \cap \mathbb{C}[e_{ij}, f_{i'j'}]$, where $\JJ_{k,n,n'}$ is the ideal of the polynomial ring $\mathbb{C}[e_{ij}, f_{i'j'},x_s,y_{s'},z_{ss'}]$. 
Such ideal, i.e., the ideal obtained by eliminating some variables, can be computed by using the theory of \textit{Gr\"obner basis}. 
For the detail, we refer the reader to \cite[Section 3]{CLO}. 
}\end{rem}

%%%%%%%%%%%%%%%%%%%%%%%%%%%%%%%%%%%%%%%%%%%%%%%%%%%%%%%%%%%%%%%%%%%%%%%%%%%%%%%%%%%%%%%%%%%%%%%%%%%%%%%%%%%%%%%%%%%%%%%%
%%%%%%%%%%%%%%%%%%%%%%%%%%%%%%%%%%%%%%%%%%%%%%%%%%%%%%%%%%%%%%%%%%%%%%%%%%%%%%%%%%%%%%%%%%%%%%%%%%%%%%%%%%%%%%%%%%%%%%%%
%%%%%%%%%%%%%%%%%%%%%%%%%%%%%%%%%%%%%%%%%%%%%%%%%%%%%%%%%%%%%%%%%%%%%%%%%%%%%%%%%%%%%%%%%%%%%%%%%%%%%%%%%%%%%%%%%%%%%%%%
\subsection{Refinement of Conjecture~\ref{conj1} using a characterization of critical graphs}\label{sec-main3}
%%%%%%%%%%%%%%%%%%%%%%%%%%%%%%%%%%%%%%%%%%%%%%%%%%%%%%%%%%%%%%%%%%%%%%%%%%%%%%%%%%%%%%%%%%%%%%%%%%%%%%%%%%%%%%%%%%%%%%%%
%%%%%%%%%%%%%%%%%%%%%%%%%%%%%%%%%%%%%%%%%%%%%%%%%%%%%%%%%%%%%%%%%%%%%%%%%%%%%%%%%%%%%%%%%%%%%%%%%%%%%%%%%%%%%%%%%%%%%%%%
%%%%%%%%%%%%%%%%%%%%%%%%%%%%%%%%%%%%%%%%%%%%%%%%%%%%%%%%%%%%%%%%%%%%%%%%%%%%%%%%%%%%%%%%%%%%%%%%%%%%%%%%%%%%%%%%%%%%%%%%

Now we consider an additional condition that
\vspace{-.5\baselineskip}
\begin{enumerate}
\setlength{\parskip}{0cm}
\setlength{\itemsep}{0cm}
\item[{\bf (V8)}]
$G$ and $H$ are $k$-critical graphs.
\end{enumerate}
\vspace{-.5\baselineskip}
Let 
$$
V'_{k,n,n'}=\{(G,H)\in V_{k,n,n'} : (G,H) \text{ satisfies (V8)}\}. 
$$
Then the following holds. 

\begin{prop}%%%%%%%%%%%%%%%%%%%%%%%%%%%%%%%%%%%%%%%%%%%%%%%%%%%%%%%%%%%%%%%%%%%%%%%%%%%%%%%%%%%%%%%%%%%%%%%%%%%%%%%%%%%%
\label{prop-main3-1}
For integers $k\geq 2$, $n\geq 1$ and $n'\geq 1$, the conditions (H1)--(H3) in Proposition~\ref{prop-main1-1} are equivalent to
\vspace*{-2.5mm}
\begin{enumerate}[{\bf (H$'$1)}]
\setlength{\parskip}{0cm}
\setlength{\itemsep}{0cm}
\item
$V'_{k,n,n'}\subseteq W_{k,n,n'}$.
\end{enumerate}
\vspace*{-2.5mm}
\end{prop}
%%%%%%%%%%%%%%%%%%%%%%%%%%%%%%%%%%%%%%%%%%%%%%%%%%%%%%%%%%%%%%%%%%%%%%%%%%%%%%%%%%%%%%%%%%%%%%%%%%%%%%%%%%%%%%%%%%%%%%%%
\proof
By Proposition~\ref{prop-main1-1}, it suffices to show that ``(H$'$1)~$\Leftrightarrow $~(H1)''.
``(H1)~$\Rightarrow$~(H$'$1)'' clearly holds.
Thus we suppose that (H$'$1) holds and show that (H1) holds.

Let $(G,H)\in V_{k,n,n'}$.
Then $G$ and $H$ satisfy (X1) and (V1)--(V7).
If $G$ and $H$ are $k$-critical, then $(G,H)\in V'_{k,n,n'}$, and so $(G,H)\in W_{k,n,n'}$ because (H$'$1) holds.
Thus, without loss of generality, we may assume that $G$ is not $k$-critical.
Since $G$ satisfies (V1), this implies that $\chi (G)\leq k-1$.
In particular, $\min\{\chi (G),\chi (H)\}\leq k-1$.
Hence $(G,H)\in W_{k,n,n'}$, which proves that (H1) holds.
\qed

Therefore, if $k$-critical graphs of order at most $n$ can be characterized, then the information for such critical graphs directly effect the system of equation corresponding to $V'_{k,n,n'}$.
By Theorem~\ref{thm-small-main}, the smallest nontrivial case for Conjecture~\ref{conj1} is $\chi (G)=\chi (H)=5$, $|V(G)|=8$ and $|V(H)|=11$.

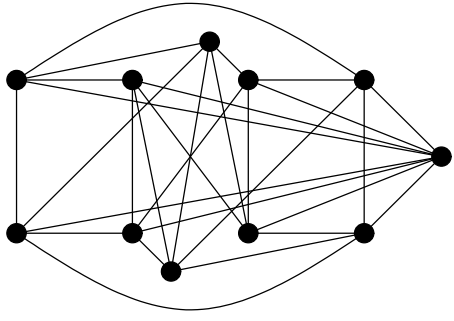
\begin{figure}[htbp]
\begin{center}
%WinTpicVersion4.32a
{\unitlength 0.1in%
\begin{picture}(23.0000,16.0000)(13.5000,-32.0000)%
% CIRCLE 2 0 0 0 Black Black  
% 4 1400 2000 1400 2050 1400 2050 1400 2050
% 
\special{sh 1.000}%
\special{ia 1400 2000 50 50 0.0000000 6.2831853}%
\special{pn 8}%
\special{ar 1400 2000 50 50 0.0000000 6.2831853}%
% CIRCLE 2 0 0 0 Black Black  
% 4 1400 2800 1400 2850 1400 2850 1400 2850
% 
\special{sh 1.000}%
\special{ia 1400 2800 50 50 0.0000000 6.2831853}%
\special{pn 8}%
\special{ar 1400 2800 50 50 0.0000000 6.2831853}%
% CIRCLE 2 0 0 0 Black Black  
% 4 2000 2800 2000 2850 2000 2850 2000 2850
% 
\special{sh 1.000}%
\special{ia 2000 2800 50 50 0.0000000 6.2831853}%
\special{pn 8}%
\special{ar 2000 2800 50 50 0.0000000 6.2831853}%
% CIRCLE 2 0 0 0 Black Black  
% 4 2000 2000 2000 2050 2000 2050 2000 2050
% 
\special{sh 1.000}%
\special{ia 2000 2000 50 50 0.0000000 6.2831853}%
\special{pn 8}%
\special{ar 2000 2000 50 50 0.0000000 6.2831853}%
% CIRCLE 2 0 0 0 Black Black  
% 4 2600 2000 2600 2050 2600 2050 2600 2050
% 
\special{sh 1.000}%
\special{ia 2600 2000 50 50 0.0000000 6.2831853}%
\special{pn 8}%
\special{ar 2600 2000 50 50 0.0000000 6.2831853}%
% CIRCLE 2 0 0 0 Black Black  
% 4 2600 2800 2600 2850 2600 2850 2600 2850
% 
\special{sh 1.000}%
\special{ia 2600 2800 50 50 0.0000000 6.2831853}%
\special{pn 8}%
\special{ar 2600 2800 50 50 0.0000000 6.2831853}%
% CIRCLE 2 0 0 0 Black Black  
% 4 3200 2800 3200 2850 3200 2850 3200 2850
% 
\special{sh 1.000}%
\special{ia 3200 2800 50 50 0.0000000 6.2831853}%
\special{pn 8}%
\special{ar 3200 2800 50 50 0.0000000 6.2831853}%
% CIRCLE 2 0 0 0 Black Black  
% 4 3200 2000 3200 2050 3200 2050 3200 2050
% 
\special{sh 1.000}%
\special{ia 3200 2000 50 50 0.0000000 6.2831853}%
\special{pn 8}%
\special{ar 3200 2000 50 50 0.0000000 6.2831853}%
% CIRCLE 2 0 0 0 Black Black  
% 4 2400 1800 2400 1850 2400 1850 2400 1850
% 
\special{sh 1.000}%
\special{ia 2400 1800 50 50 0.0000000 6.2831853}%
\special{pn 8}%
\special{ar 2400 1800 50 50 0.0000000 6.2831853}%
% CIRCLE 2 0 0 0 Black Black  
% 4 2200 3000 2200 3050 2200 3050 2200 3050
% 
\special{sh 1.000}%
\special{ia 2200 3000 50 50 0.0000000 6.2831853}%
\special{pn 8}%
\special{ar 2200 3000 50 50 0.0000000 6.2831853}%
% CIRCLE 2 0 0 0 Black Black  
% 4 3600 2400 3600 2450 3600 2450 3600 2450
% 
\special{sh 1.000}%
\special{ia 3600 2400 50 50 0.0000000 6.2831853}%
\special{pn 8}%
\special{ar 3600 2400 50 50 0.0000000 6.2831853}%
% LINE 2 0 3 0 Black White  
% 16 1400 2000 1400 2800 1400 2800 2000 2800 2000 2800 2600 2000 2600 2000 3200 2000 3200 2000 3200 2800 3200 2800 2600 2800 2600 2800 2000 2000 2000 2000 1400 2000
% 
\special{pn 8}%
\special{pa 1400 2000}%
\special{pa 1400 2800}%
\special{fp}%
\special{pa 1400 2800}%
\special{pa 2000 2800}%
\special{fp}%
\special{pa 2000 2800}%
\special{pa 2600 2000}%
\special{fp}%
\special{pa 2600 2000}%
\special{pa 3200 2000}%
\special{fp}%
\special{pa 3200 2000}%
\special{pa 3200 2800}%
\special{fp}%
\special{pa 3200 2800}%
\special{pa 2600 2800}%
\special{fp}%
\special{pa 2600 2800}%
\special{pa 2000 2000}%
\special{fp}%
\special{pa 2000 2000}%
\special{pa 1400 2000}%
\special{fp}%
% SPLINE 2 0 3 0 Black White  
% 4 1400 2000 2300 1600 3200 2000 3200 2000
% 
\special{pn 8}%
\special{pa 1400 2000}%
\special{pa 1429 1981}%
\special{pa 1458 1961}%
\special{pa 1488 1942}%
\special{pa 1517 1922}%
\special{pa 1546 1903}%
\special{pa 1575 1885}%
\special{pa 1605 1866}%
\special{pa 1692 1812}%
\special{pa 1722 1795}%
\special{pa 1751 1778}%
\special{pa 1809 1746}%
\special{pa 1839 1731}%
\special{pa 1868 1716}%
\special{pa 1897 1702}%
\special{pa 1926 1689}%
\special{pa 1956 1677}%
\special{pa 1985 1665}%
\special{pa 2014 1654}%
\special{pa 2043 1644}%
\special{pa 2073 1635}%
\special{pa 2102 1627}%
\special{pa 2131 1620}%
\special{pa 2160 1614}%
\special{pa 2190 1609}%
\special{pa 2219 1605}%
\special{pa 2248 1602}%
\special{pa 2277 1600}%
\special{pa 2307 1600}%
\special{pa 2336 1601}%
\special{pa 2365 1603}%
\special{pa 2394 1606}%
\special{pa 2423 1611}%
\special{pa 2453 1616}%
\special{pa 2511 1630}%
\special{pa 2540 1639}%
\special{pa 2570 1648}%
\special{pa 2628 1670}%
\special{pa 2657 1682}%
\special{pa 2687 1695}%
\special{pa 2716 1708}%
\special{pa 2745 1723}%
\special{pa 2774 1737}%
\special{pa 2804 1753}%
\special{pa 2862 1785}%
\special{pa 2891 1802}%
\special{pa 2921 1820}%
\special{pa 3008 1874}%
\special{pa 3038 1893}%
\special{pa 3125 1950}%
\special{pa 3155 1970}%
\special{pa 3184 1989}%
\special{pa 3200 2000}%
\special{fp}%
% SPLINE 2 0 3 0 Black White  
% 4 1400 2800 2300 3200 3200 2800 3200 2800
% 
\special{pn 8}%
\special{pa 1400 2800}%
\special{pa 1429 2819}%
\special{pa 1458 2839}%
\special{pa 1488 2858}%
\special{pa 1517 2878}%
\special{pa 1546 2897}%
\special{pa 1575 2915}%
\special{pa 1605 2934}%
\special{pa 1692 2988}%
\special{pa 1722 3005}%
\special{pa 1751 3022}%
\special{pa 1809 3054}%
\special{pa 1839 3069}%
\special{pa 1868 3084}%
\special{pa 1897 3098}%
\special{pa 1926 3111}%
\special{pa 1956 3123}%
\special{pa 1985 3135}%
\special{pa 2014 3146}%
\special{pa 2043 3156}%
\special{pa 2073 3165}%
\special{pa 2102 3173}%
\special{pa 2131 3180}%
\special{pa 2160 3186}%
\special{pa 2190 3191}%
\special{pa 2219 3195}%
\special{pa 2248 3198}%
\special{pa 2277 3200}%
\special{pa 2307 3200}%
\special{pa 2336 3199}%
\special{pa 2365 3197}%
\special{pa 2394 3194}%
\special{pa 2423 3189}%
\special{pa 2453 3184}%
\special{pa 2511 3170}%
\special{pa 2540 3161}%
\special{pa 2570 3152}%
\special{pa 2628 3130}%
\special{pa 2657 3118}%
\special{pa 2687 3105}%
\special{pa 2716 3092}%
\special{pa 2745 3077}%
\special{pa 2774 3063}%
\special{pa 2804 3047}%
\special{pa 2862 3015}%
\special{pa 2891 2998}%
\special{pa 2921 2980}%
\special{pa 3008 2926}%
\special{pa 3038 2907}%
\special{pa 3125 2850}%
\special{pa 3155 2830}%
\special{pa 3184 2811}%
\special{pa 3200 2800}%
\special{fp}%
% LINE 2 0 3 0 Black White  
% 8 2400 1800 1400 2800 1400 2000 2400 1800 2400 1800 2600 2000 2400 1800 2600 2800
% 
\special{pn 8}%
\special{pa 2400 1800}%
\special{pa 1400 2800}%
\special{fp}%
\special{pa 1400 2000}%
\special{pa 2400 1800}%
\special{fp}%
\special{pa 2400 1800}%
\special{pa 2600 2000}%
\special{fp}%
\special{pa 2400 1800}%
\special{pa 2600 2800}%
\special{fp}%
% LINE 2 0 3 0 Black White  
% 8 2200 3000 2000 2800 2200 3000 2000 2000 2200 3000 3200 2800 2200 3000 3200 2000
% 
\special{pn 8}%
\special{pa 2200 3000}%
\special{pa 2000 2800}%
\special{fp}%
\special{pa 2200 3000}%
\special{pa 2000 2000}%
\special{fp}%
\special{pa 2200 3000}%
\special{pa 3200 2800}%
\special{fp}%
\special{pa 2200 3000}%
\special{pa 3200 2000}%
\special{fp}%
% LINE 2 0 3 0 Black White  
% 16 3600 2400 3200 2000 3600 2400 2600 2000 2000 2000 3600 2400 3600 2400 1400 2000 1400 2800 3600 2400 3600 2400 2000 2800 2600 2800 3600 2400 3600 2400 3200 2800
% 
\special{pn 8}%
\special{pa 3600 2400}%
\special{pa 3200 2000}%
\special{fp}%
\special{pa 3600 2400}%
\special{pa 2600 2000}%
\special{fp}%
\special{pa 2000 2000}%
\special{pa 3600 2400}%
\special{fp}%
\special{pa 3600 2400}%
\special{pa 1400 2000}%
\special{fp}%
\special{pa 1400 2800}%
\special{pa 3600 2400}%
\special{fp}%
\special{pa 3600 2400}%
\special{pa 2000 2800}%
\special{fp}%
\special{pa 2600 2800}%
\special{pa 3600 2400}%
\special{fp}%
\special{pa 3600 2400}%
\special{pa 3200 2800}%
\special{fp}%
% LINE 2 0 3 0 Black White  
% 4 2000 2000 2000 2800 2600 2800 2600 2000
% 
\special{pn 8}%
\special{pa 2000 2000}%
\special{pa 2000 2800}%
\special{fp}%
\special{pa 2600 2800}%
\special{pa 2600 2000}%
\special{fp}%
% LINE 2 0 3 0 Black White  
% 2 2400 1800 2200 3000
% 
\special{pn 8}%
\special{pa 2400 1800}%
\special{pa 2200 3000}%
\special{fp}%
\end{picture}}%
\caption{A $K_{4}$-free $5$-critical graph $H^{*}$ of order $11$.}
\label{f-JR}
\end{center}
\end{figure}

Jensen and Royle \cite{JR} also claimed that there exist 56 $K_{4}$-free graphs with chromatic number $5$. 
Although they are not always $5$-critical, we expect that there are a lot of $5$-critical ones among them. %and so it seems to be difficult to list all such graphs.
Now, we focus on the graph $H^{*}$ depicted in Figure~\ref{f-JR}, which is one of $K_4$ free $5$-critical graph of order $11$, 
and demonstrate a partial solution of the problem that $V'_{5,8,11}\subseteq W_{5,8,11}$. 

We translate the situations into ideal type formulas. All the ideals below are in the polynomial ring 
$$\mathbb{C}[e_{ij},f_{i'j'},z_{ss'} : 1 \leq i < j \leq 8, 1 \leq i' < j' \leq 11, 1 \leq s \leq 8, 1 \leq s' \leq 11]$$
of $28+55+88=171$ variables. Let
\begin{align}\tag{W1($G=K_1+H_0, k=5$)}
\begin{split}
E &= (e_{12}-1,e_{13}-1,e_{14}-1,e_{25}-1,e_{26}-1,e_{35}-1,e_{37}-1,e_{46}-1,e_{47}-1,e_{56}-1,e_{57}-1,e_{67}-1)\\
&+ (e_{i8}-1:1\leq i\leq 7)+(e_{ij}:\mbox{otherwise}) 
\end{split}
\end{align}
and 
\begin{align}\tag{W1($H=H^{*}, k=5$)}
\begin{split}
E' &= (f_{12}-1,f_{13}-1,f_{17}-1,f_{24}-1,f_{28}-1,f_{34}-1,f_{36}-1,f_{45}-1,f_{56}-1,f_{57}-1,f_{68}-1,f_{78}-1)\\
&+ (f_{19}-1,f_{29}-1,f_{59}-1,f_{69}-1,f_{3,10}-1,f_{4,10}-1,f_{7,10}-1,f_{8,10}-1,f_{9,10}-1)\\
&+ (f_{i,11}-1:1\leq i\leq 8)+(f_{ij}:\mbox{otherwise}).
\end{split}
\end{align}
Let 
$$
\LL=E+E'+Z_{5,8,11}+J_{5,8,11}.
$$
To show the equality 
\begin{align}\label{mainaim-sec5}
\LL=\mathbb{C}[e_{ij},f_{i'j'},z_{ii'}], 
\end{align}
which implies that $\LL \cap \mathbb{C}[e_{ij},f_{i'j'}] = \mathbb{C}[e_{ij},f_{i'j'}]$, 
gives a partial solution of the problem that $V'_{5,8,11}\subseteq W_{5,8,11}$. 
We will focus on \eqref{mainaim-sec5} and related problems in Section~\ref{sec-nume}.
Remark that to get a complete solution of $V'_{5,8,11}\subseteq W_{5,8,11}$, 
it suffices to prove a similar inclusion problem for all $K_{4}$-free $5$-critical graphs $H$ except for $H^{*}$.

%%%%%%%%%%%%%%%%%%%%%%%%%%%%%%%%%%%%%%%%%%%%%%%%%%%%%%%%%%%%%%%%%%%%%%%%%%%%%%%%%%%%%%%%%%%%%%%%%%%%%%%%%%%%%%%%%%%%%%%%
%%%%%%%%%%%%%%%%%%%%%%%%%%%%%%%%%%%%%%%%%%%%%%%%%%%%%%%%%%%%%%%%%%%%%%%%%%%%%%%%%%%%%%%%%%%%%%%%%%%%%%%%%%%%%%%%%%%%%%%%
%%%%%%%%%%%%%%%%%%%%%%%%%%%%%%%%%%%%%%%%%%%%%%%%%%%%%%%%%%%%%%%%%%%%%%%%%%%%%%%%%%%%%%%%%%%%%%%%%%%%%%%%%%%%%%%%%%%%%%%%
\section{Computational experiment}\label{sec-nume}
%%%%%%%%%%%%%%%%%%%%%%%%%%%%%%%%%%%%%%%%%%%%%%%%%%%%%%%%%%%%%%%%%%%%%%%%%%%%%%%%%%%%%%%%%%%%%%%%%%%%%%%%%%%%%%%%%%%%%%%%
%%%%%%%%%%%%%%%%%%%%%%%%%%%%%%%%%%%%%%%%%%%%%%%%%%%%%%%%%%%%%%%%%%%%%%%%%%%%%%%%%%%%%%%%%%%%%%%%%%%%%%%%%%%%%%%%%%%%%%%%
%%%%%%%%%%%%%%%%%%%%%%%%%%%%%%%%%%%%%%%%%%%%%%%%%%%%%%%%%%%%%%%%%%%%%%%%%%%%%%%%%%%%%%%%%%%%%%%%%%%%%%%%%%%%%%%%%%%%%%%%
For confirming that Conjecture~\ref{conj1} is true or finding a counterexample in the case the graphs are small, 
we implement Theorem~\ref{thm3-1} and other related functions by an open source general computer algebra system Risa/Asir \cite{ASIR}. 
All source codes of our programming are put at the webpage \cite{URL}. 
More precisely, we implement the computations whether the following inclusion or the equality are true or not: 
\begin{itemize}
\setlength{\parskip}{0cm}
\setlength{\itemsep}{0cm}
\item[(a)] $\tilde{\JJ}_{k,n,n'}\subseteq \tilde{\II}_{k,n,n'}$ for given $n,n'$ and $k$; %(see Theorem~\ref{thm3-1}); 
\item[(b)] $\LL=\mathbb{C}[e_{ij},f_{i'j'},z_{ii'}]$; %(see \eqref{mainaim-sec5}). 
\end{itemize}
Note that \cite{URL} contains many other functions related to our problem, some of which will be explained below. 
Theorem~\ref{thm3-1} says that the inclusion (a) is equivalent to that Conjecture~\ref{conj1} is true. 
On the other hand, the discussions developed in Section~\ref{sec-main3} say that 
the confirmations of the equality (b) implies the search of the smallest non-trivial unknown case of Conjecture~\ref{conj1}.

For example, we can perform the computations in Mac or Linux OS as follows: 
\begin{Ex}{\em 
  Let $k = 3, n = 5, n' = 5$. 
  On Risa/Asir running with terminal, we will check the inclusion (a), i.e., we will check the condition Theorem~\ref{thm3-1} (ii) as follows: 
  \begin{verbatim}
    [1895] load("**certan path**/Hedetniemi.rr")$
    [1949] K = 3$
    [1950] N = 5$
    [1951] N' = 5$
    [1952] hedetniemi.theorem_4_4(N, N', K);
    ### Theorem 4.4 (ii): k = 3, n = 5, n' = 5
    (omitted)
    ### True: k = 3, n = 5, n' = 5
  \end{verbatim}
Similarly, we can check the equality (b) by \verb|hedetniemi.section_4_3()|. 
}\end{Ex}

We performed the above computational experiments (a). 
All computations have been performed in Ubuntu OS equipped with 64 GB memory, Intel Xeon(R) W-2135, CPU 3.7 GHz. %CPU *, ** GHz. 
The following table shows the times took for each experiment. 
As the tables show, it took a huge time for checking even trivial cases of Conjecture~\ref{conj1}.

\begin{center}
\begin{tabular}{c|c|c}
                     & Time        &Result \\ \hline\hline
    $k=3, n=4, n'=4$ & $2$ seconds & True \\ \hline
    $k=3, n=4, n'=5$ & $58$ seconds & True \\ \hline
    $k=3, n=5, n'=5$ & $4409$ seconds $\fallingdotseq 73$ minutes & True \\ \hline
    $k=3, n=5, n'=6$ & $241116$ seconds $\fallingdotseq 67$ hours  & True \\ \hline
    $k=3, n=6, n'=6$ & more than two weeks & Still running \\ \hline
    $k=4, n=4, n'=4$ & $167$ seconds & True \\ \hline
    $k=4, n=4, n'=5$ & $480818$ seconds $\fallingdotseq 133$ hours & True \\ \hline
    $k=4, n=5, n'=5$ & more than one month & Still running \\
\end{tabular}
\end{center}

Unfortunately, the computation (b) did not stop even after one month. 
We have to upgrade the machine performance or devise the algorithm or the theoretical part in order to push the boundary of the computable cases 
(more concretely, to complete the coputation of \verb|hedetniemi.section_4_3()|).

Instead, we implemented the following experimental computations: 
\begin{itemize}
\item[(c-1)] we replace $E$ (resp. $E'$) with the ideal corresponding to $C_{2m+1}$ (resp. $C_{2m'+1}$), and $\LL=E+E'+Z_{3,2m+1,2m'+1}+J_{3,2m+1,2m'+1}$; 
\item[(c-2)] we replace both $E$ and $E'$ with the ideal corresponding to the graph $H_0$, and $\LL=E+E'+Z_{4,7,7}+J_{4,7,7}$; 
\item[(c-3)] we replace $E$ (resp. $E'$) with the ideal corresponding to the graph $H_0$ (resp. the graph depicted in Figure~\ref{fig-Gr}), and $\LL=E+E'+Z_{4,7,11}+J_{4,7,11}$; 
\end{itemize}
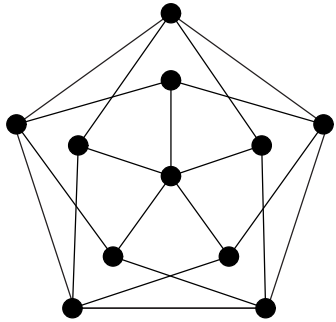
\begin{figure}[htbp]
\begin{center}
%WinTpicVersion4.32a
{\unitlength 0.1in%
\begin{picture}(16.9000,16.4000)(42.5000,-36.5000)%
% POLYGON 2 0 3 0 Black White  
% 6 5600 3600 4600 3600 4291 2649 5100 2061 5909 2649 5600 3600
% 
\special{pn 8}%
\special{pa 5600 3600}%
\special{pa 4600 3600}%
\special{pa 4291 2649}%
\special{pa 5100 2061}%
\special{pa 5909 2649}%
\special{pa 5600 3600}%
\special{pa 4600 3600}%
\special{fp}%
% CIRCLE 2 0 0 0 Black Black  
% 4 5100 2910 5100 2960 5100 2960 5100 2960
% 
\special{sh 1.000}%
\special{ia 5100 2910 50 50 0.0000000 6.2831853}%
\special{pn 8}%
\special{ar 5100 2910 50 50 0.0000000 6.2831853}%
% CIRCLE 2 0 0 0 Black Black  
% 4 5100 2410 5100 2460 5100 2460 5100 2460
% 
\special{sh 1.000}%
\special{ia 5100 2410 50 50 0.0000000 6.2831853}%
\special{pn 8}%
\special{ar 5100 2410 50 50 0.0000000 6.2831853}%
% CIRCLE 2 0 0 0 Black Black  
% 4 5100 2060 5100 2110 5100 2110 5100 2110
% 
\special{sh 1.000}%
\special{ia 5100 2060 50 50 0.0000000 6.2831853}%
\special{pn 8}%
\special{ar 5100 2060 50 50 0.0000000 6.2831853}%
% CIRCLE 2 0 0 0 Black Black  
% 4 5890 2640 5890 2690 5890 2690 5890 2690
% 
\special{sh 1.000}%
\special{ia 5890 2640 50 50 0.0000000 6.2831853}%
\special{pn 8}%
\special{ar 5890 2640 50 50 0.0000000 6.2831853}%
% CIRCLE 2 0 0 0 Black Black  
% 4 4300 2640 4300 2690 4300 2690 4300 2690
% 
\special{sh 1.000}%
\special{ia 4300 2640 50 50 0.0000000 6.2831853}%
\special{pn 8}%
\special{ar 4300 2640 50 50 0.0000000 6.2831853}%
% CIRCLE 2 0 0 0 Black Black  
% 4 4590 3600 4590 3650 4590 3650 4590 3650
% 
\special{sh 1.000}%
\special{ia 4590 3600 50 50 0.0000000 6.2831853}%
\special{pn 8}%
\special{ar 4590 3600 50 50 0.0000000 6.2831853}%
% CIRCLE 2 0 0 0 Black Black  
% 4 5590 3600 5590 3650 5590 3650 5590 3650
% 
\special{sh 1.000}%
\special{ia 5590 3600 50 50 0.0000000 6.2831853}%
\special{pn 8}%
\special{ar 5590 3600 50 50 0.0000000 6.2831853}%
% CIRCLE 2 0 0 0 Black Black  
% 4 5570 2750 5570 2800 5570 2800 5570 2800
% 
\special{sh 1.000}%
\special{ia 5570 2750 50 50 0.0000000 6.2831853}%
\special{pn 8}%
\special{ar 5570 2750 50 50 0.0000000 6.2831853}%
% CIRCLE 2 0 0 0 Black Black  
% 4 4620 2750 4620 2800 4620 2800 4620 2800
% 
\special{sh 1.000}%
\special{ia 4620 2750 50 50 0.0000000 6.2831853}%
\special{pn 8}%
\special{ar 4620 2750 50 50 0.0000000 6.2831853}%
% CIRCLE 2 0 0 0 Black Black  
% 4 4800 3330 4800 3380 4800 3380 4800 3380
% 
\special{sh 1.000}%
\special{ia 4800 3330 50 50 0.0000000 6.2831853}%
\special{pn 8}%
\special{ar 4800 3330 50 50 0.0000000 6.2831853}%
% CIRCLE 2 0 0 0 Black Black  
% 4 5400 3330 5400 3380 5400 3380 5400 3380
% 
\special{sh 1.000}%
\special{ia 5400 3330 50 50 0.0000000 6.2831853}%
\special{pn 8}%
\special{ar 5400 3330 50 50 0.0000000 6.2831853}%
% LINE 2 0 3 0 Black White  
% 10 5400 3330 5100 2910 5100 2910 5570 2750 5100 2910 5100 2410 5100 2910 4620 2750 5100 2910 4800 3330
% 
\special{pn 8}%
\special{pa 5400 3330}%
\special{pa 5100 2910}%
\special{fp}%
\special{pa 5100 2910}%
\special{pa 5570 2750}%
\special{fp}%
\special{pa 5100 2910}%
\special{pa 5100 2410}%
\special{fp}%
\special{pa 5100 2910}%
\special{pa 4620 2750}%
\special{fp}%
\special{pa 5100 2910}%
\special{pa 4800 3330}%
\special{fp}%
% LINE 2 0 3 0 Black White  
% 20 4800 3330 4300 2640 4300 2640 5100 2410 5100 2410 5890 2640 5890 2640 5400 3330 5400 3330 4590 3600 4590 3600 4620 2750 4620 2750 5100 2060 5100 2060 5570 2750 5570 2750 5590 3600 5590 3600 4800 3330
% 
\special{pn 8}%
\special{pa 4800 3330}%
\special{pa 4300 2640}%
\special{fp}%
\special{pa 4300 2640}%
\special{pa 5100 2410}%
\special{fp}%
\special{pa 5100 2410}%
\special{pa 5890 2640}%
\special{fp}%
\special{pa 5890 2640}%
\special{pa 5400 3330}%
\special{fp}%
\special{pa 5400 3330}%
\special{pa 4590 3600}%
\special{fp}%
\special{pa 4590 3600}%
\special{pa 4620 2750}%
\special{fp}%
\special{pa 4620 2750}%
\special{pa 5100 2060}%
\special{fp}%
\special{pa 5100 2060}%
\special{pa 5570 2750}%
\special{fp}%
\special{pa 5570 2750}%
\special{pa 5590 3600}%
\special{fp}%
\special{pa 5590 3600}%
\special{pa 4800 3330}%
\special{fp}%
\end{picture}}%
\caption{A triangle-free $4$-critical graph of order $11$.}
\label{fig-Gr}
\end{center}
\end{figure}
We note that $C_{2m+1}$ is $3$-critical, $H_0$ is $4$-critical of order $11$. 
We also note that we already know theoritcally that the equality \eqref{mainaim-sec5} holds for any case (c-1)--(c-3). 
Those computations stopped within some minutes or hours as shown below: 
\begin{center}
\begin{tabular}{c|c|c|c}
$E$      &$E'$     &Time        &Result \\ \hline\hline
%$C_{11}$ &$C_{11}$ &134 seconds &True \\ \hline
%$C_{11}$ &$C_{13}$ &459 seconds &True \\ \hline 
$C_{13}$ &$C_{13}$ &1270 seconds $\fallingdotseq 21$ minutes &True \\ \hline 
$C_{13}$ &$C_{15}$ &3855 seconds $\fallingdotseq 64$ minutes &True \\ \hline 
$C_{15}$ &$C_{15}$ &9890 seconds $\fallingdotseq 164$ minutes &True \\ \hline 
$C_{15}$ &$C_{17}$ &27743 seconds $\fallingdotseq 7.7$ hours &True \\ \hline 
$C_{17}$ &$C_{17}$ &63158 seconds $\fallingdotseq 17.5$ hours &True \\
\end{tabular}
\quad
\begin{tabular}{c|c|c|c}
$E$   &$E'$  &Time &Result \\ \hline\hline
$H_0$ &$H_0$ &1 second &True \\ \hline
$H_0$ &Figure~\ref{fig-Gr} &18 seconds &True \\
\end{tabular}
\end{center}

%%%%%%%%%%%%%%%%%%%%%%%%%%%%%%%%%%%%%%%%%%%%%%%%%%%%%%%%%%%%%%%%%%%%%%%%%%%%%%%%%%%%%%%%%%%%%%%%%%%%%%%%%%%%%%%%%%%%%%%%
%%%%%%%%%%%%%%%%%%%%%%%%%%%%%%%%%%%%%%%%%%%%%%%%%%%%%%%%%%%%%%%%%%%%%%%%%%%%%%%%%%%%%%%%%%%%%%%%%%%%%%%%%%%%%%%%%%%%%%%%
%%%%%%%%%%%%%%%%%%%%%%%%%%%%%%%%%%%%%%%%%%%%%%%%%%%%%%%%%%%%%%%%%%%%%%%%%%%%%%%%%%%%%%%%%%%%%%%%%%%%%%%%%%%%%%%%%%%%%%%%
\section{Concluding remark}\label{sec-conclu}
%%%%%%%%%%%%%%%%%%%%%%%%%%%%%%%%%%%%%%%%%%%%%%%%%%%%%%%%%%%%%%%%%%%%%%%%%%%%%%%%%%%%%%%%%%%%%%%%%%%%%%%%%%%%%%%%%%%%%%%%
%%%%%%%%%%%%%%%%%%%%%%%%%%%%%%%%%%%%%%%%%%%%%%%%%%%%%%%%%%%%%%%%%%%%%%%%%%%%%%%%%%%%%%%%%%%%%%%%%%%%%%%%%%%%%%%%%%%%%%%%
%%%%%%%%%%%%%%%%%%%%%%%%%%%%%%%%%%%%%%%%%%%%%%%%%%%%%%%%%%%%%%%%%%%%%%%%%%%%%%%%%%%%%%%%%%%%%%%%%%%%%%%%%%%%%%%%%%%%%%%%

In this paper, we presented a reduction of Conjecture~\ref{conj1} using the inclusion of the ideals of a polynomial ring (Theorem~\ref{thm3-1}).
Since our reduction strongly depends on the structure of critical graphs as we verified in Subsection~\ref{sec-main3}, the advance of the research of the criticality directly gives favorable effects on Conjecture~\ref{conj1}.

We remark that Shitov~\cite{Sh} used the existence of graphs with large fractional chromatic number and large girth, and so his counterexamples implicitly depend on so-called probabilistic method.
In particular, it seems to be difficult to give their specific constructions. 
Since our main result (Theorem~\ref{thm3-1}) gives a reduction for each case, we expect that it offers not only a new approach for Conjecture~\ref{conj1} (in small chromatic number case) but the smallest specific counterexample of Conjecture~\ref{conj1}.

Furthermore, every Shitov's counterexample contains a large clique.
Hence the following weaker conjecture than Conjecture~\ref{conj1} is naturally posed.

\begin{conj}%%%%%%%%%%%%%%%%%%%%%%%%%%%%%%%%%%%%%%%%%%%%%%%%%%%%%%%%%%%%%%%%%%%%%%%%%%%%%%%%%%%%%%%%%%%%%%%%%%%%%%%%%%%%
\label{conj2}
Let $G$ and $H$ be triangle-free graphs.
Then $\chi (G\times H)=\min\{\chi (G),\chi (H)\}$.
\end{conj}
%%%%%%%%%%%%%%%%%%%%%%%%%%%%%%%%%%%%%%%%%%%%%%%%%%%%%%%%%%%%%%%%%%%%%%%%%%%%%%%%%%%%%%%%%%%%%%%%%%%%%%%%%%%%%%%%%%%%%%%%

Conjecture~\ref{conj2} is still interesting because the chromatic number of graphs with large girth has deeply studied in graph theory.
Note that the triangle-freeness of a labeled graph on $[n]$ is corresponding the following condition:
\begin{align}
e_{ij}e_{jl}e_{il}=0~~~(1\leq i<j<l\leq n),\label{eq-conclu-1}
\end{align}
where $e_{ij}$, $e_{jl}$ and $e_{il}$ are in Subsection~\ref{sec-main2}.
Since every subgraph of a triangle-free graph is also triangle-free, the criticality argument in Subsection~\ref{sec-main3} can work if we consider the triangle-free graphs.
Consequently, our reduction can be applied to Conjecture~\ref{conj2}.

\end{document}